[1994/12/01]% LaTeX date must December 1994 or later
\documentclass[reqno,13pt]{amsart}
\usepackage{bbm}
\def\1{\mathbbm{1}}

\usepackage{amsmath,amsthm}
\usepackage[mathscr]{euscript}
\usepackage{graphicx}
\usepackage{amsmath,amssymb}
\usepackage{circuitikz}

\usepackage{pstricks}
\usepackage{pst-node}

\usepackage{tikz}
\usetikzlibrary{arrows,chains,matrix,positioning,scopes}
\tikzstyle{element}=[rectangle,draw,fill=white, line width=1pt]
\tikzstyle{terminal}=[circle,draw, scale=0.3, line width=1pt,red]
\tikzstyle{fleche}=[->,>=stealth', very thick]
\tikzstyle{fleche1}=[->,>=stealth', very thick, red]
\usepackage{placeins}
\definecolor{darkgoldenrod4}{rgb}{0.55,0.4,0.55}
\definecolor{maroon4}{rgb}{0.55,0.11,0.38}
\definecolor{indianred}{rgb}{0.8,0.36,0.36}
\definecolor{purple1}{rgb}{0.61,0.19,1}
\definecolor{goldenrod1}{rgb}{1,0.76,0.15}
\definecolor{indianred3}{rgb}{0.8,0.33,0.33}
\definecolor{red4}{rgb}{0.55,0,0}
\definecolor{darkslategray}{rgb}{0.18,0.31,0.31}
\definecolor{firebrick}{rgb}{0.7,0.13,0.13}
\definecolor{slateblue3}{rgb}{0.41,0.35,0.8}
\definecolor{mediumorchid4}{rgb}{0.48,0.22,0.55}
\definecolor{thistle4}{rgb}{0.55,0.48,0.55}
\definecolor{rltred}{rgb}{0.75,0,0}
\definecolor{rltgreen}{rgb}{0,0.5,0}
\definecolor{oneblue}{rgb}{0,0,0.75}
\definecolor{marron}{rgb}{0.64,0.16,0.16}
\definecolor{forestgreen}{rgb}{0.13,0.54,0.13}
\definecolor{purple}{rgb}{0.62,0.12,0.94}
\definecolor{dockerblue}{rgb}{0.11,0.56,0.98}
\definecolor{freeblue}{rgb}{0.25,0.41,0.88}
\definecolor{myblue}{rgb}{0,0.2,0.4}

\def\R{{\mathbb{R}}}

\def\al{\alpha}

\def\ga{\gamma}

\def\la{\lambda}

\def\t{\tau}

\def\B{\mathbb{B}}

\def\calL{{\mathcal{L}}}

\def\calT{{\mathcal{T}}}

\def\calA{{\mathcal{A}}}

\def\calC{{\mathcal{C}}}

\def\R{\mathbb R}

\def\P{\mathbb P}
\def\D{\mathbb D}
\def\S{\mathbb S}
\def\C{\mathbb C}

\def\N{\mathbb N}

\def\F{\mathbb F}

\def\L{\mathbb L}

\def\T{\mathbb T}

\newtheorem{definition}{Definition}[section]

\newtheorem{lemma}{Lemma}[section]
\newtheorem{proposition}{Proposition}[section]
\newtheorem{corollary}{Corollary}[section]
\newtheorem{theorem}{Theorem}[section]
\theoremstyle{remark}
\newtheorem{remark}{Remark}[section]
\newtheorem{example}{Example}[section]
\newtheorem{mainassumptions}{Assumption}[section]

\title{Well-posedness and stability of boundary delay equations} %Semigroup generation properties of positive operator

\thanks{This work was supported by National Key Research and Development Program of China (2024YFE0214000), National Natural Science Foundation of China (Grant No. 62173308), Natural Science Foundation of Zhejiang Province of China (Grant No. LRG25F030002), Zhejiang Province Leading Geese Plan (Grant No. 2025C01056).}% At most 5 thanks
\author{Yassine El Gantouh} \address{Yassine El Gantouh, School of Mathematical Sciences, Zhejiang Normal University, Jinhua, Zhejiang 321004; and School of Mathematics, Southwest Jiaotong University, Chengdu, Sichuan 611756, P.R. China}\email{elgantouhyassine@gmail.com}

\author{Yang Liu} \address{Yang Liu, School of Mathematical Sciences, Zhejiang Normal University; China-Mozambique Belt and Road Joint Laboratory on Smart Agriculture, Zhejiang Normal University; and Hangzhou School of Automation, Zhejiang Normal University, Jinhua, Zhejiang 321004, P.R. China}\email{liuyang@zjnu.edu.cn}

 %93D20=Asymptotic stability in control theory, 93C43=delay observation\control

\begin{document}
	\maketitle
	
	\renewcommand{\sectionmark}[1]{}
	%In this paper, we investigate the well-poseness and stability of a wide class of linear abstract boundary delay problems. 
	\begin{abstract}
    In this paper, we introduce the notion of boundary delay equations, establishing a unified framework for analyzing linear time-invariant systems with pure time-delayed boundary conditions. We establish mild sufficient conditions for the existence, uniqueness, and positivity of solutions. Furthermore, we derive spectral criteria for exponential stability. The conditions on the perturbation generalize well-known criteria for the generation of domain perturbations of positive semigroup generators. As an application, we present necessary and sufficient conditions for the exponential stability of positive hyperbolic systems with time-delayed boundary conditions.
	\end{abstract}

	\section{Introduction}
    Time-delay systems arise naturally in engineering and applied sciences due to transmission delays in information, signals, or material flows.  These systems are mathematically modeled using delay differential equations, a class of functional differential equations that account for the system's dependence on its past states \cite{Hale}. The standard analytical approach maps the finite-dimensional dynamics into infinite-dimensional function spaces, such as the space of square-integrable Lebesgue functions or continuous functions. This framework allows for a unified treatment of time-delay systems \cite{CZ,fridman2014introduction,Hale}. 
    \vspace{.01cm}
    
    Over the past few decades, partial differential equations (PDEs) with state delays (SDs) have attracted considerable attention. These systems are valuable for modeling realistic and complex behaviors \cite{Wu}. Significant progress has been made in understanding their mathematical properties, including stability analysis, bifurcation phenomena, and long-term dynamics (see, e.g., \cite{BP1,BP2,BFR,Datko,EN,JGH,MVo,Webb,Wu,Xu} and references therein). Much of this research has focused on PDEs with internal SDs, where semigroup methods have been effective. However, boundary SDs pose greater challenges. The key difficulty arises because boundary SDs introduce unbounded operators into the abstract system formulation (see, e.g., Remark \ref{Observation} below). This complicates the direct application of semigroup methods and other standard analytical techniques (see, e.g., \cite{AABM,BFR1,BBH,DLP,HMR}).
    \vspace{.01cm}
    
    In this work, we introduce and analyze a novel class of delay differential equations called \emph{boundary delay equations}. This system provides an effective framework for modeling linear time-invariant (LTI) systems with time delays appearing either in boundary conditions or in feedback laws used for boundary stabilization of PDEs. We refer to \cite{AABM,Boujijane2024,HK,Zhang} for concrete applications that can potentially be studied in light of the approach we develop in our work.
    \vspace{.01cm}
    
    In particular, our study exploits the theory of delay equations. More precisely, we use product spaces and operator matrices to reformulate the boundary delay equation as a Cauchy problem. The resulting problem is governed by an operator matrix with a perturbed domain in the sense of \cite{Gr} (see also \cite{HMR}). By applying a version of the Weiss-Staffans perturbation theorem developed for positive systems, we prove that the underlying differential operator generates a positive C$_0$-semigroup. Moreover, we establish necessary and sufficient conditions for exponential stability and derive explicit spectral criteria for two important classes of boundary delay equations. 
    %\vspace{.01cm}
    
    The main novelties of this work are:
    \begin{enumerate}
    	\item We present a version of the Weiss-Staffans perturbation theorem for positive systems. The characteristic feature of this modified theorem is that it relaxes the regularity requirements for the open-loop system while preserving its positivity under static feedback via a linear bounded positive operator.
    	\item We prove the well-posedness and positivity of solutions of boundary delay equations under rather weak regularity assumptions on the system operators. 
    	\item We provide spectral criteria for the exponential stability of boundary delay equations.
    \end{enumerate}
    \vspace{.01cm}
    
    An important motivation for this work is to provide spectral criteria for exponential stability of linear positive hyperbolic systems subjected to time-delayed proportional boundary feedback laws \cite{Zhang}. While such feedback mechanisms are well-established for stabilizing hyperbolic systems \cite{BaCo}, our contribution extends these results by incorporating time delays into proportional boundary feedback laws. We note that \cite{BBH} investigated the exponential stability of linear positive hyperbolic systems with time-delayed boundary conditions bounded with respect to the state variable \cite[Sect. 5]{BBH}. Consequently, their proposed boundary coupling corresponds to a full-state feedback control, which requires measuring the state over the entire spatial domain. This requirement can be impractical in real-world applications (see, e.g., \cite{Zhang}). 
    \vspace{.01cm}
    
    This article is organized as follows. First, relevant notation and the problem statement are given in Section \ref{Sec:1}. In Section \ref{Sec:2}, we present a version of the Weiss-Staffans perturbation theorem for positive systems. This section begins with preliminary results on positive admissible control and observation operators in Subsections \ref{Sec:2.1} and \ref{Sec:2.2}. We then state and prove the Weiss-Staffans perturbation theorem in Subsection \ref{Sec:2.3}. Section \ref{Sec:3} contains the main results of this work; Theorem \ref{Main1-23} establishes sufficient conditions for the well-posedness and positivity of boundary delay equations, while Theorem \ref{Main1-24} and Corollaries \ref{Cor2} and \ref{Cor3} provide spectral criteria for exponential stability.  In Section \ref{Sec:4}, we apply our theoretical results to derive necessary and sufficient conditions for the exponential stability of positive hyperbolic systems with time-delayed boundary conditions. Finally, in Appendix \ref{Appendix}, we present generation and stability theorems for domain perturbations of positive semigroup generators.
	
	\section{Notation and Problem Statement}\label{Sec:1}
	\subsection{Notation and Terminology}\label{Sec:1.1}
   	Let $\mathbb{C},\mathbb{R},\mathbb{R}_+,\mathbb{N}$ denote the sets of complex numbers, real numbers, positive real numbers, and natural numbers, respectively. Let $E$ and $F$ be Banach spaces. We denote by $\calL(E,F)$ the space of bounded linear operators from $E$ to $F$, and we write $\calL(E):=\calL(E,E)$. For an operator $P\in \mathcal{L}(E,F)$ and a subspace $Z$ of $E$, we denote by $P\vert_{Z}$ the restriction of $P$ on $Z$. The range of the operator $P$ is denoted by Range$(P)$. For a real $p\ge 1$, we denote by $L^p(\R_+;E)$ the space of all measurable $p$-integrable functions $f:\R_+\to E$. For $a,b\in \R$ such that $b\ge a$, we denote by $\calC(a,b;E)$ the space of all continuous functions from $[a,b]$ to $E$, while by $W^{1,p}(a,b;E)$ we denote the Sobolev space of absolutely continuous functions $f: [a,b] \to E$ such that the derivative $\partial_x f$ is a $p$-integrable function.  The shorthand notations $\partial_t:=\tfrac{\partial}{\partial t}$ and $\partial_\theta:=\tfrac{\partial}{\partial \theta}$ are utilized for the partial derivatives.
    \vspace{.01cm}
    
    A \emph{Banach lattice} $E:=(E,\le)$ is a partially ordered Banach space in which every pair of elements $x,y$ of $E$ have a supremum $x \vee y$, and for all $x,y,z\in E$ and $\alpha\ge 0$:
   \begin{eqnarray*}
   	\left\lbrace
   	\begin{array}{lll}
   		x\leq y &\implies (x+z\le y+z \; {\rm and} \; \alpha x\leq \alpha y ),\\
   		\vert x\vert \leq \vert y\vert & \implies \Vert x\Vert \leq \Vert y\Vert,
   	\end{array}
   	\right.
   \end{eqnarray*}
   where $\vert x\vert =x\vee (-x)$ for all $x\in E$. An element $x\in E$ is called positive if $x\ge 0$, and the set of all positive elements forms the positive cone $E_+$ (a notation also used for general ordered Banach spaces). An operator $ P\in \mathcal{L}(E,F)$ is called positive if it maps the positive cone $E_+$ into $F_+$, where $F$ is a Banach lattice. The set of all such positive operators is denoted by $\mathcal{L}_+(E,F)$. This notation is also used when $E$ and $F$ are order Banach spaces. For further details, see \cite{CHARALAMBOS} or \cite{Schaf}.
   
   \subsection{Positive  C$_0$-semigroups and Extrapolation Spaces}\label{Sec:1.2}
   Let $A$ be the generator of a C$_0$-semigroup $\mathbb{T}:=(T(t))_{t\geq0}$ on the Banach space $E$, with growth bound $\omega_0(A)$ and spectral bound $s(A)$. Define $E_{-1}^{A}$ as the completion of $E$ under the norm $\|x\|_{-1}^{A}:=\|R(\la,A)x\|$, where $\la\in \rho(A)$ (the resolvent set of $A$) is fixed and $R(\la,A):=(\la I_E-A)^{-1}$.  The space $E_{-1}^{A}$ is independent of $\lambda$, as different choices yield equivalent norms; see, e.g., \cite[Sect. II.5]{EN}. The semigroup $\mathbb{T}$ extends to a C$_0$-semigroup $\mathbb{T}_{-1} := (T_{-1}(t))_{t \geq 0}$ on $E_{-1}^{A}$, with generator $A_{-1}$, an extension of $A$ with domain $E$. 
   \vspace{.1cm}
   
   Now, let $E$ be a Banach lattice. The semigroup $\T$ is called positive if $T(t)x\ge 0$ for all $x\in E_+$ and $t\ge 0$. The operator $A$ is called resolvent positive if there exists $\omega\in \mathbb{R}$ such that $(\omega,\infty)\subset \rho (A) $ and $R(\la,A)\geq 0$ for all $\la > \omega$. While every generator of a positive C$_0$-semigroup is resolvent positive, the converse is not true in general. Counterexamples and further discussion can be found in, e.g.,  \cite[Sect. 3]{Arendt}.  For an overview of the theory of positive C$_0$-semigroups, we refer to the monographs \cite{BFR,Nagel}. 
   \vspace{.1cm}
   
   We note that in our setting, $E$ is a Banach lattice with a natural positive cone $E_+$. However, it is not immediately clear how to extend this concept of positivity to the extrapolation space $E_{-1}^A$. Since $A$ generates a positive C$_0$-semigroup on $E$, we follow \cite[Def. 2.1]{BJVW} and define the positive cone in $E_{-1}^A$, denoted by $(E_{-1}^A)_+$,  as the closure of $E_+$ in the norm $\|\cdot\|_{-1}^{A}$. This definition ensures that $E_+\subset (E_{-1}^A)_+$. Furthermore, if $E$ is a real Banach lattice, then by \cite[Prop. 2.3]{BJVW} we have
   \begin{align*}
   	E_+=E\cap (E_{-1}^{A})_+.
   \end{align*}
   Further details can be found in \cite[Sect. 2.2]{SGPS} and \cite{SGF}.

	\subsection{Problem Statement}\label{Sec:1.3}
	In this subsection, we introduce the concept of boundary delay equations that can be used to describe the behavior of a wide class of dynamical systems. Let $X$ and $V$ be Banach spaces equipped with norms $\Vert \cdot\Vert_X$ and $\Vert \cdot\Vert_V$. For $r\in \R_+$ and $1\le p<\infty$, let $Y_p$ be defined by $Y_p:=L^p(-r,0;V)$ and be
	endowed with the norm
	\begin{align*}
		\Vert \varphi\Vert_{Y_p}:=\int_{-r}^{0}\Vert \varphi(\theta)\Vert^p_{V}d\theta.
	\end{align*}
	On this space we introduce the operator $L:W^{1,p}(-r,0;V)\to V,$ called delay operator, by
	\begin{align}\label{dealy-operator}
		L(\varphi):=\int_{-r}^{0} d\eta(\theta) \varphi(\theta),
	\end{align}
	where $\varphi\in W^{1,p}(-r,0;V)$ and $\eta:[-r,0]\to \calL(V)$ is an operator-valued function that satisfies the following assumptions.
	\begin{mainassumptions}\label{Assp22}
		\begin{itemize}
			\item[(i)] $\eta$ is of bounded variation, i.e.,  
			\begin{align*}
				{\rm Var}(\eta;-r,0):=\sup\left \{ \sum_{k=1}^{m}\Vert \eta(\theta_k)-\eta(\theta_{k-1})\Vert_{\calL(V)},
				\theta_1=-r<\ldots<\theta_m=0,\; m\in \N\right\}<+\infty;
			\end{align*}
			\item[(ii)] $\eta$ is continuous from the left, i.e.,
			\begin{align*}
				\lim_{\theta \to \theta_0^{-}}\eta(\theta)=\eta(\theta_0^{-}).
			\end{align*}
		\end{itemize}
	\end{mainassumptions}
	\begin{remark}\label{further1}
		In the above, (i) implies that the Stieltjes integral $\int_{-r}^{0} d\eta(\theta) \varphi(\theta)$ exists for every $\varphi\in \calC(-r,0;V)$ and is defined as the limit of Riemann-Stieltjes sums:
		\begin{align*}
			\sum_{k=1}^{m}\left( \eta(\theta_k)-\eta(\theta_{k-1})\right)\varphi(a_k) \; {\rm as}\; \max_k\vert\theta_k-\theta_{k-1}\vert\mapsto 0,
		\end{align*}
		where $a_k\in [\theta_{k-1},\theta_{k}]$. While (ii) yields that 
		\begin{align*}
			\lim_{\alpha\to 0}{\rm Var}(\eta;-\alpha,0)=0.
		\end{align*}
		Additionally, $L(\varphi)$ is well-defined for all $\varphi\in W^{1,p}(-r,0;V)$, since by the Sobolev embedding theorem we know that $W^{1,p}(-r,0;V)$ is densely and continuously embedded in $\calC(-r,0;V)$; see, e.g., \cite{MVo} for more information.
	\end{remark}
	
	Let $A_m: D(A_m) \subseteq X \to X$ be a closed, densely defined linear operator, and let $G, P: D(A_m) \to V$ be continuous linear boundary operators. With this notation, a boundary delay problem is described by the following LTI system
	\begin{align}\label{Eq1}
		\begin{cases}
			\partial_t z(t)=A_mz(t),&t> 0,\\  
			Gz(t)= L(x_t), & t> 0,\\
			x(t)=Pz(t),& t>0,
		\end{cases}
	\end{align}
	with initial states $z(0)=f\in X$ and $x_0=\varphi\in Y_p$. 
	Here, $ z(t) $ is the state and $x_t$ is the so-called history segment defined as 
	\begin{align*}
		[-r,0]\ni \theta\mapsto x_t(\theta)=x(t+\theta)\in V.
	\end{align*}
	
	\begin{remark}
		We note that the function $x_t:[-r,0]\to V$ describes the history segment of the value of the state $z(t)$ at the boundary. When $P$ is a bounded operator, we recover the delay systems studied in \cite[Ex. 5.2]{HMR} and \cite[Sect. 5]{BBH}. Thus, the delay differential equation \eqref{Eq1} generalizes boundary delay equations by incorporating modifications of the history function through an unbounded boundary operator. This framework is particularly relevant for the stability analysis of LTI systems with pure time delays, whether in boundary conditions or in the feedback laws employed for boundary stabilization of PDEs.
	\end{remark}
	
	In the sequel, we make the following standard assumptions.
	\begin{mainassumptions}\label{Assp11}
		\begin{itemize}
			\item[(i)] The restriction $A:=(A_m)\vert_{\ker(G)}$ generates a C$_0$-semigroup $\T:=(T(t))_{t\ge 0}$ on $X$;
			\item[(ii)] $G$ is onto.  
		\end{itemize}
	\end{mainassumptions}
	Under these assumptions, the domain of $A_m$ can be decomposed and related to $A$ as follows
	\begin{eqnarray*}
		D(A_m)=\ker (G) \oplus \ker ( \la I_X- A_{m} ),\qquad \la\in \rho(A).
	\end{eqnarray*}
	Furthermore, for any $\lambda\in\rho(A)$, the Dirichlet operator associated with the pair $(A_m,G)$,  defined by
	\begin{eqnarray*}%\label{Dirichlet}
		\D_\lambda:=\left(G_{|\ker(\lambda I_X-A_m)}\right)^{-1}: V\to \ker(\lambda I_X-A_m)\subseteq X,
	\end{eqnarray*} 
	exists and is bounded (see \cite{Gr,Salam} for more details). One can then verify that the operator
	\begin{eqnarray*}
		\B:=(\lambda I_X-A_{-1})\D_\lambda\in\mathcal{L}(V,X_{-1}^A),
	\end{eqnarray*}
	is independent of $\lambda$, Range$(\B)\cap X=\{0\}$, and that 
	\begin{eqnarray}\label{representation}
		(A_m-A_{-1})|_{D(A_m)}=\B G.
	\end{eqnarray}
	\begin{remark}\label{Observation}
		Using \eqref{representation}, the boundary delay equation \eqref{Eq1} can be reformulated as
		\begin{align}\label{Delay-equation}
			\begin{cases}
				\partial_t z(t)=A_{-1}z(t)+\B L(x_t),&t> 0,\cr  
				x(t)=Pz(t),& t>0.
			\end{cases}
		\end{align}
		For the case $P=I_X$ and $\B=I_X$, the well-posedness and asymptotic behavior of \eqref{Delay-equation} were established in \cite{BP1,BP2} using the Miyadera-Voigt perturbation theorem. When $\B$ is unbounded on $X$ and $A$ generates an analytic C$_0$-semigroup, \cite[Sect. 3.4]{BP2} proved well-posedness via the same perturbation technique. This result was extended in \cite{JGH} for discrete delay operators $L$ without requiring the semigroup to be analytic. More recently, \cite{BFR1} addressed well-posedness and positivity properties of solutions for the boundary delay equation \eqref{Eq1} with $P=I_X$ and $G=\tilde{G}-N$, where $\tilde{G}$ satisfies Assumption \ref{Assp11} and $N$ is unbounded operator that is not necessarily closed or closable. In their analysis, the authors combined perturbation theory of semigroups with feedback theory of linear infinite-dimensional systems (see \cite{HMR}) to establish the well-posedness and positivity of solutions.
	\end{remark}
	
	The typical approach to the study of delay differential equations is to embed the underlying dynamical process into certain product space. To this end, we define
	\begin{eqnarray*}
		\mathscr{X}_p :=X\times Y_{p},\qquad \left\|(f, \varphi)^{\top}\right\|:=\|f\|_{X}+\|\varphi\|_{Y_p},
	\end{eqnarray*}
	and the new state
	\begin{align*}%\label{new-state}
		\zeta(t)=(z(t), x_{t})^{\top},\qquad t\ge 0.
	\end{align*}
	We point out here that the function $x_t $ is the solution of the transport equation
	\begin{align}\label{shift}
		\begin{cases}
			\partial_t x_t= \partial_\theta x_t,& t> 0,\\
			x_t(0)= x(t), &t> 0.
		\end{cases}
	\end{align} 
	With this notation, the abstract delay equation \eqref{Eq1} can be rewritten as the following Cauchy problem
	\begin{align}\label{Cauchy-problem}
		\begin{cases}
			\partial_t\zeta(t)=\mathfrak{A}\zeta(t),& t>0, \\ \zeta(0)=(f,\varphi)^{\top},&
		\end{cases}
	\end{align}
	where $\mathfrak{A}:D(A_m)\times  W^{1,p}(-r,0;V) \to  \mathscr{X}_p$ is the operator matrix defined by
	\begin{align}\label{cauchy-operator}
	\mathfrak{A}:={\rm diag}(A_m ,\partial_\theta), \qquad D(\mathfrak{A}):=\left\{ f\in D(A_m),\varphi\in  W^{1,p}(-r,0;V):\ Gf=L(\varphi),\;\varphi(0)=Pf\right\}.
	\end{align}
	Our goal is then to show that the operator $(\mathfrak{A},D(\mathfrak{A}))$ generates a C$_0$-semigroup on $\mathscr{X}_p$ and that there exists $\alpha>0$ such that
	\begin{align*}
		\left\Vert (z(t),x_t)^\top\right\Vert_{\mathscr{X}_p}\le \Upsilon_{f,\varphi} \,e^{-\alpha t}\left\Vert (f,\varphi)^\top\right\Vert_{ \mathscr{X}_p},
	\end{align*}
	for all $(f,\varphi)^\top\in D(\mathfrak{A})$, $t\ge 0$, and a constant $\Upsilon_{f,\varphi}>0$. 
	
	\section{Preliminary results on linear time-invariant infinite-dimensional positive systems}\label{Sec:2}
		Our proof of the main results is based on the feedback theory of linear infinite-dimensional positive systems. In the following subsections, we present auxiliary results essential to our analysis. Throughout this section, $\mathsf{X}$ and $\mathsf{U}$ are Banach lattices (we use the same notation $\Vert \cdot\Vert$ to specify the corresponding norms), and $\mathsf{A}$ is the generator of a positive C$_0$-semigroup $\mathsf{T}:=(\mathsf{T}(t))_{t\ge 0}$ on $\mathsf{X}$.
	
	\subsection{Admissible positive observation operators}\label{Sec:2.1}
	Consider the observed linear system
	\begin{align}\label{observed-system}
		\begin{cases}
			\partial_t z(t)=\mathsf{A} z(t), & t>0, \quad z(0)=x,\cr y(t)=\mathsf{C} z(t), & t> 0,
		\end{cases}
	\end{align}
	where $\mathsf{C}\in \calL(D(\mathsf{A}),\mathsf{U})$. For $x\in D(\mathsf{A})$, the solution $z(t)=\mathsf{T}(t)x$ and the observation function $y(t)=\mathsf{C}\mathsf{T}(t)x$ are well-defined for all $t\ge 0$, since $D(\mathsf{A})$ is invariant under $\mathsf{T}$. However, defining $y(\cdot,x)$ for arbitrary $x \in \mathsf{X}$ is not straightforward. The system \eqref{observed-system} is called well-posed if the output $y$ can be extended to a function $y\in L^p_{loc}(\R_+,\mathsf{U})$ satisfying
	\begin{align*}
		\|y(\cdot,x)\|_{L^p(0,\al;\mathsf{U})}\le \ga\|x\|,\qquad \forall x\in \mathsf{X},\,\alpha>0,
	\end{align*}
	where $\ga:=\ga(\al)>0$ is a constant. 
	
	We select the following definition; see, e.g., \cite[Sect. 4.3]{TW}.
	\begin{definition}\label{S2.D2}
		For a real number $p\ge 1$, we say that $\mathsf{C}$ is an $L^p$-admissible observation operator for $ \mathsf{A}$  (or $(\mathsf{C},\mathsf{A})$ is $L^p$-admissible for short) if for some (hence for all) $\alpha >0,$ there exists a constant $\gamma:=\gamma(\alpha)>0$ such that
		\begin{align}\label{S2.7}
			\int_{0}^{\alpha} \Vert \mathsf{C}\mathsf{T}(t)x\Vert^{p}dt\leq \gamma^p\Vert x\Vert^{p},\qquad \forall x\in D(\mathsf{A}). 
		\end{align}
		%	If $\lim_{\alpha\mapsto 0 }\gamma(\alpha)=0 $, then $(\mathsf{C},\mathsf{A})$ is called zero-class $ L^p $-admissible.
	\end{definition}
	\begin{remark}\label{Remark-output-map}
		We  note that if $(\mathsf{C},\mathsf{A})$ is $L^p$-admissible, then the  operator
		\begin{align*}
			\Psi: D(\mathsf{A})\to L^p_{loc}(\R_+,\mathsf{U}),\quad x\mapsto \Psi x:= \mathsf{C}\mathsf{T}(\cdot)x,
		\end{align*}
		called the output-map of $(\mathsf{C},\mathsf{A})$, is well defined and has a bounded extension to $\tilde{\Psi}:\mathsf{X}\to L^p(0,\al;\mathsf{U})$ for each $\al>0$. We associate with $\mathsf{C}$ the following operator
		\begin{align}\label{cgama}
			\begin{split}
				D(\mathsf{C}_{\Lambda})&:=\{x\in \mathsf{X}:\lim_{\lambda\to+\infty} \mathsf{C}\lambda R(\lambda,\mathsf{A})x\;{\rm exists}\},\cr
				\mathsf{C}_{\Lambda}x& :=\lim_{\lambda \to+ \infty} \mathsf{C}\la  R(\la ,\mathsf{A})x,\quad x\in D(\mathsf{C}_{\Lambda}),
			\end{split}
		\end{align}
		called the Yosida extension of $\mathsf{C}$ w.r.t $\mathsf{A}$. The operator introduced above makes possible to give a simple pointwise interpretation of the extended  output function $y$ in terms of $\mathsf{C}$. More precisely, we have 
		\begin{align*}
			y(t,x)= \mathsf{C}_\Lambda \mathsf{T}(t)x, 
		\end{align*}
		for all $x\in \mathsf{X}$ and a.e. $t\ge 0$. For further details; see, e.g., \cite[Sect. 4.4]{Staf}.
	\end{remark}
	The following result demonstrates that the admissibility of a positive observation operator for positive semigroups is entirely characterized by its behavior on $D_+(A)$ \cite[Lem. 3.1]{El} (see also \cite{SGPS}).
	\begin{proposition}\label{positive-admis-C}
		Let $ \mathsf{C}\in\mathcal{L}_+(D(\mathsf{A}),\mathsf{U})$ and $p\ge 1$ be a real number. If the estimate \eqref{S2.7} holds for all $ x\in D_+(\mathsf{A})$ and for some $ \alpha > 0 $, then $(\mathsf{C},\mathsf{A})$ is $L^p$-admissible.
	\end{proposition}
	\begin{remark}
		With respect to the above result one can reduce the definition of admissibility of positive observation operators  (more generally, the difference of two positive operators) for positive semigroups to 
		\begin{align*}%\label{S2.7}
			\int_{0}^{\alpha} \Vert \mathsf{C}\mathsf{T}(t)x\Vert^{p}dt\leq \gamma^{p}\Vert x\Vert^{p},\qquad \forall\,  x\in D_+(\mathsf{A}),
		\end{align*}
		for some $\al>0$ and a constant $\gamma:=\ga(\al)>0$. In this case, we say that $(\mathsf{C},\mathsf{A})$ is positive $L^p$-admissible. 
		
		We note that if $\mathsf{C}$ is positive, then, by the closedness of $\mathsf{U}_+$, the Yosida extension $\mathsf{C}_{\Lambda}$ of $\mathsf{C}$ w.r.t $\mathsf{A}$ is also positive. Thus,  for $(\mathsf{C},\mathsf{A})$ positive $L^p$-admissible, one obtains that $0\le y(t,x)= \mathsf{C}_{\Lambda}\mathsf{T}(t)x $ for all $ x\in \mathsf{X}_+ $ and a.e. $ t\geq 0 $.
	\end{remark}
	
	\subsection{Admissible positive control operators}\label{Sec:2.2}
	Consider the controlled equation
	\begin{align}\label{diff-equ}
		\begin{cases}
			\partial_t z(t)=\mathsf{A}_{-1} z(t) +\mathsf{B}v(t),& t>0,\\
			z(0)=x, &
		\end{cases}
	\end{align}
	where $\mathsf{A}_{-1}:\mathsf{X}\to \mathsf{X}_{-1}^\mathsf{A}$ is the extension of $\mathsf{A}$ to $\mathsf{X}$ and $ \mathsf{B}\in \mathcal{L}(\mathsf{U},\mathsf{X}_{-1}^\mathsf{A}) $. A mild solution of \eqref{diff-equ} is given by 
	\begin{align}\label{S2.2}
		z(t,x,v)=\mathsf{T}(t)x+\int_{0}^{t}\mathsf{T}_{-1}(t-s)\mathsf{B}v(s)ds,
	\end{align}
	for all $t\ge 0$, $ x\in \mathsf{X} $, and $v\in  L^{p}_{loc}(\R_+;\mathsf{U}) $. The integral in \eqref{S2.2} is calculated in $ \mathsf{X}_{-1}^\mathsf{A} $. However for practice purpose, we typically seek continuous $\mathsf{X}$-valued functions $z$. This requires identifying a class of control operators $\mathsf{B}$ for which the state of \eqref{diff-equ} remains in $\mathsf{X}$. In particular, we select the following definition (see, e.g., \cite{WC}).
	\begin{definition}\label{S2.D1}
		For a real number $p\ge 1$, we say that $ \mathsf{B}\in \mathcal{L}(\mathsf{U},\mathsf{X}_{-1}^\mathsf{A}) $ is an $ L^p $-admissible control operator for $\mathsf{A}$ (or $(\mathsf{A},\mathsf{B})$ is $L^p$-admissible for short) if for some (and hence for all) $t>0$ and all $ v\in L^{p}(\R_+;\mathsf{U}) $
		\begin{align*}%\label{S2.2'}
			\Phi_{t}v:=\int_{0}^{t}\mathsf{T}_{-1}(t-s)\mathsf{B}v(s)ds\in \mathsf{X}.
		\end{align*}
	\end{definition}
	We note that by the closed graph theorem if $(\mathsf{A},\mathsf{B})$ is $L^p$-admissible, then for every $t>0$ there exists a constant $\kappa:=\kappa(t)>0$ such that
	\begin{align}\label{input-map}
		\Vert \Phi_{t}v\Vert \leq \kappa\Vert v\Vert_{ L^{p}(0,\t;\mathsf{U})},\quad \forall v\in L^{p}(\R_+;\mathsf{U}) .
	\end{align}
	As a matter of fact, the solutions of \eqref{diff-equ} are continuous $\mathsf{X}$-valued functions of $t$ for any $x\in \mathsf{X}$ and any $v\in L^p_{loc}(\R_+,\mathsf{U})$ with $p\in [1,+\infty)$, cf. \cite[Prop. 2.3]{WC}.  Moreover, by \cite[Prop. 2.1]{El}, the following holds.
	\begin{proposition}\label{state-trajec}
		Let $ \mathsf{B}\in \mathcal{L}_+(\mathsf{U},\mathsf{X}_{-1}^\mathsf{A})$ and $p\in [1,+\infty)$. If for some $t>0$
		\begin{align*}
			\Phi_t v\in  \mathsf{X}_+,\qquad \forall v\in L^p_+(\R_+;\mathsf{U}),
		\end{align*}
		then $(\mathsf{A},\mathsf{B})$ is called positive $L^p$-admissible. If this is the case, then the state $z$ of \eqref{diff-equ} satisfy $z\in \calC(\R_+,\mathsf{X})$ and
		\begin{eqnarray*}
			\mathsf{X}_+\ni z(t)= \mathsf{T}(t)x+ \Phi_tv ,\; \forall t\ge 0, x\in \mathsf{X}_+, v\in L^{p}_+(\R_+;\mathsf{U}). 
		\end{eqnarray*}
	\end{proposition}
	
	\subsection{A Weiss-Staffans perturbation theorem for positive systems}\label{Sec:2.3}   
	Consider the following input-output system
	\begin{eqnarray}\label{input-ouptut}
		\begin{cases}
			\partial_t z(t) =\mathsf{A}_{-1} z(t)+\mathsf{B}v(t),& t> 0,\quad z(0)=x,\\
			y(t) = \Gamma z(t), & t> 0,
		\end{cases}
	\end{eqnarray}
	where $\mathsf{B}\in \calL(\mathsf{U},\mathsf{X}_{-1}^\mathsf{A})$ and $ \Gamma\in \calL(\mathsf{Z},\mathsf{U})$ with
	\begin{align*}
		\mathsf{Z}=D(\mathsf{A})+R(\la,\mathsf{A}_{-1})\mathsf{B}\mathsf{U}, \qquad \la\in \rho(\mathsf{A}).
	\end{align*}
	It is not difficult to check that $\mathsf{Z}$ endowed with the norm
	\begin{eqnarray*}
		\Vert z\Vert_{\mathsf{Z}}^2=\inf\left\{\Vert x\Vert^2+\Vert v\Vert^2: x\in D(\mathsf{A}),v\in \mathsf{U},\ z=x+R(\la,\mathsf{A}_{-1}) \mathsf{B}v \right\},
	\end{eqnarray*}
	is a Banach space such that $D(\mathsf{A})\subseteq \mathsf{Z}\subseteq \mathsf{X}$ with continuous embeddings. We note that $\mathsf{Z}$ is independent of the choice of $\la$, since different choices of $\la$ lead to equivalent norms on $\mathsf{Z}$ due to the resolvent equation. 
	
	In what follows, we assume that $(\mathsf{A},\mathsf{B})$ is positive $L^p$-admissible. By Proposition \ref{state-trajec}, we have $ z(t,x,v)\in \mathsf{X}$ for all $t\ge 0$, $x\in \mathsf{X}$, and $v\in L^{p}(\R_+;\mathsf{U})$. However, the observation function $ y(\cdot) $ is not a priori well-defined, since $\Gamma\in \calL(\mathsf{Z},\mathsf{U})$ and $z$ may not belong to $\mathsf{Z}$. To address this issue, we introduce the space
	\begin{align*}
		W^{1,p}_{0}(0,\t;\mathsf{U}):=\left\{v\in W^{1,p}(0,\t;\mathsf{U}):v(0)=0\right\},
	\end{align*}
	for a fixed $\tau>0$. Assuming without loss of generality that $0\in \rho(\mathsf{A})$ and using integration by parts, it follows that $\Phi_t v\in \mathsf{Z} $ for all $v\in W^{1,p}_{0}(0,\tau;\mathsf{U})$ and $t\in [0,\t]$. This allows us to define the following operator for such $v$:
	\begin{align*}
		(\mathbb{F}v)(t):=\Gamma\Phi_t v, \qquad t\in [0,\t].
	\end{align*}
	Using this framework, we obtain the output representation
	\begin{eqnarray}\label{output-fct}
		y(\cdot,x,v)= \Psi x+\F v, \quad \forall x\in D(\mathsf{A}),\, v\in W^{1,p}_{0}(0,\t;\mathsf{U}),
	\end{eqnarray}
	where $\Psi$ is the output-map of $(\mathsf{C},\mathsf{A})$ with $\mathsf{C}:=\Gamma\vert_{ D(\mathsf{A})}$.
	\vspace{.1cm}
	
	Based on this construction, we define the well-posedness of system \eqref{input-ouptut} as follows.
	\begin{definition}\label{Definition-Lp-well-posed}%\textcolor{freeblue}{
			The input-output system \eqref{input-ouptut} is called $L^p$-well-posed (for $1\le p<\infty$) if, for every $ \tau > 0 $, initial state $x\in D(\mathsf{A})$, and input $v\in W^{1,p}_{0}(0,\tau;\mathsf{U})$, there exists a constant $ c_\tau>0 $ (independent of $x$ and $v$) such that the following estimate holds for all solutions $(z,y)$ of \eqref{input-ouptut}:
			\begin{eqnarray}\label{Sigma}
				\Vert z(\tau)\Vert^{p}_{\mathsf{X}}+\Vert y\Vert^{p}_{ L^{p}(0,\tau;\mathsf{U})}\leq c_\tau\left( \Vert x\Vert^{p}_{\mathsf{X}}+\Vert v\Vert^{p}_{ L^{p}(0,\tau;\mathsf{U})}\right).
			\end{eqnarray}
		\end{definition}
		The well-posedness of system \eqref{input-ouptut} involves establishing conditions on the operators $\mathsf{A}$, $\mathsf{B}$, and $\Gamma$ that ensure continuous dependence of both the output signal and final state on the initial state and input signal. This continuity enables the extension of the solution to arbitrary initial states in $\mathsf{X}$ and input signals in $L^p_{loc}(\R_+;\mathsf{U})$;  with the corresponding output signals then belonging to $L^p_{loc}(\R_+;\mathsf{U})$.
		\vspace{.1cm}
		
		The following result yields a characterization of the well-poseness of \eqref{input-ouptut} when the operators involved are positive.
		\begin{proposition}\label{ABC-sigma}
			Let $\mathsf{B}\in \calL_+(\mathsf{U},\mathsf{X}_{-1,\mathsf{A}})$ and $\Gamma\in \calL_+(\mathsf{Z},\mathsf{U})$. Then, the input-output system \eqref{input-ouptut} is $L^p$-well-posed (with $1\le p<\infty$) if and only if 
			\begin{itemize}
				\item[(\emph{i})] $(\mathsf{A},\mathsf{B}) $ is positive $L^p$-admissible;
				\item[(\emph{ii})] $(\mathsf{C},\mathsf{A})$ is positive $L^p$-admissible;
				\item[(\emph{iii})] For $ \tau >0 $, there exists $ \beta :=\beta(\tau)>0 $ such that 
				\begin{eqnarray*}\label{wellposed-estimate}
					\Vert \F v\Vert_{L^{p}(0,\tau;\mathsf{U})}\leq \beta \Vert v\Vert_{L^{p}(0,\tau;\mathsf{U})},\quad \forall v\in W^{1,p}_{0,+}(0,\t;\mathsf{U}).
				\end{eqnarray*}
			\end{itemize}
		\end{proposition}
		\begin{proof}
			Necessity. This follows directly from Definition \ref{Definition-Lp-well-posed}.
			
			Sufficiency. To establish sufficiency, we need to verify that the state $z$ and the output $y$ of \eqref{input-ouptut} satisfy the estimate \eqref{Sigma}. In fact, let $ M\geq 1 $ and $ w\geq \omega_0(\mathsf{A})$ such that $ \Vert \mathsf{T}(t)\Vert \leq Me^{wt} $ for all $ t\geq 0 $. Using item (\emph{i}), it follows from Proposition \ref{state-trajec} that
			\begin{align*}
				\|z(\t)\|^p_\mathsf{X}\le M_1\left(\|x\|_\mathsf{X}^p+\|v\|_{L^p(0,\t;\mathsf{U})}^p\right), 
			\end{align*}  
			for all $\t\ge 0$ and $(x,v)\in \mathsf{X}\times L^p(0,\tau;\mathsf{U})$, where $M_1:=M_1(\tau)>0$ depends on $\kappa,M,w$. On the other hand, we know from Proposition \ref{positive-admis-C} that item (\emph{ii}) implies that the output map $\Psi$ is extendable to $\tilde{\Psi}\in \mathcal{L}(\mathsf{X},L^{p}(0,\t;\mathsf{U}))$ for each $\tau>0$. Moreover, item (\emph{iii}) combined with the density of $W^{1,p}_{0,+}(0,\t;\mathsf{U})$ in $L^p_+(0,\t;\mathsf{U})$ ensures that $ \mathbb{F}$ extends uniquely to a bounded operator (still denoted by $\F$) from $L^{p}(0,\t;\mathsf{U})$ to itself for each $\t> 0$. Thus, by \eqref{output-fct}, one obtains 
			$$\Vert y\Vert^{p}_{ L^{p}(0,\t;\mathsf{U})}\leq M_2\left( \Vert x\Vert^{p}_{\mathsf{X}}+\Vert v\Vert^{p}_{ L^{p}(0,\t;\mathsf{U})}\right),$$ 
			for all $(x,v)\in \mathsf{X}\times L^p(0,\tau;\mathsf{U})$, where the constant $M_2:=M_2(\tau)>0$ depends on $\gamma$ and $\beta$. This completes the proof. 
		\end{proof}
		
		\begin{definition}\label{Def.triple}
			Consider the setting of Proposition \ref{ABC-sigma}. We say that $(\mathsf{A},\mathsf{B},\mathsf{C})$ is a positive $L^p$-well-posed triple on ($\mathsf{U},\mathsf{X},\mathsf{U}$) if the conditions (\emph{i})-(\emph{iii}) in Proposition \ref{ABC-sigma} are satisfied. 
		\end{definition}
		
		If $(\mathsf{A},\mathsf{B},\mathsf{C})$ is a positive $L^p$-well-posed triple, then the operator $ \mathbb{F}$ extends to  $\F\in \calL(L^{p}(0,\t;\mathsf{U}))$ for any $\t> 0$,  and the corresponding output function $y$ of the system \eqref{input-ouptut} satisfies 
		\begin{eqnarray*}
			y(\cdot,x,v)= \Psi x+\F v, \quad \forall x\in \mathsf{X},\, v\in L^p_{loc}(\R_+;\mathsf{U}).
		\end{eqnarray*}
		Furthermore, there exists a unique analytic $\mathcal{L}(\mathsf{U})$-valued function $\mathbf{H}$, called the transfer function of $(\mathsf{A},\mathsf{B},\mathsf{C})$, which satisfies (cf. \cite[Cor. 3.5]{StaffansWeiss2002})
		\begin{align*}
			\mathbf{H}(\lambda)=\Gamma R(\lambda,\mathsf{A}_{-1})\mathsf{B},\qquad \forall\ \text{Re}\lambda>\omega_0(\mathsf{A}).
		\end{align*}
		Note that $\Gamma$ is an extension of $C$, which is generally not unique. Consequently, we do not distinguish notationally between two transfer functions for  $(\mathsf{A},\mathsf{B},\mathsf{C})$ that differ only by a fixed operator, and both are analytic. For further details; see \cite{StaffansWeiss2002} or \cite[Sect. 4.6]{Staf}.
		
		In order to give a complete representation of the extended output function $y$, we introduce the following subclass of $L^p$-well-posed linear systems \cite{WR} (see also \cite[Sect. 5.6]{Staf}).
		\begin{definition}\label{regurality}
			Let $(\mathsf{A},\mathsf{B},\mathsf{C})$ be a positive $L^p$-well-posed triple with $1\le p< \infty$. 
			\begin{itemize}
				\item We say that $(\mathsf{A},\mathsf{B},\mathsf{C})$ is a positive $L^p$-well-posed weakly regular triple with zero feedthrough if 
				\begin{eqnarray*}
					\lim_{\R\ni \la\to +\infty }\langle\mathbf{H}(\la)v,u^*\rangle=0,\quad \forall v\in \mathsf{U},\,u^*\in \mathsf{U}',
				\end{eqnarray*}
				where $\mathsf{U}'$ is the topological dual of $\mathsf{U}$.
				\item We say that $(\mathsf{A},\mathsf{B},\mathsf{C})$ is a positive $L^p$-well-posed strongly regular triple with zero feedthrough if 
				\begin{eqnarray*}
					\lim_{\R\ni \la\to +\infty }\Vert\mathbf{H}(\la)v\Vert=0, \qquad \forall v\in \mathsf{U}.
				\end{eqnarray*}
				\item We say that $(\mathsf{A},\mathsf{B},\mathsf{C})$ is a positive $L^p$-well-posed uniformly regular triple with zero feedthrough if
				\begin{eqnarray*}
					\lim_{\R\ni \la\to +\infty }\Vert\mathbf{H}(\la)\Vert_{\calL(\mathsf{U})}=0.
				\end{eqnarray*} 
			\end{itemize}
		\end{definition}
		
		If, for instance, $(\mathsf{A},\mathsf{B},\mathsf{C}) $ is a positive $L^p$-well-posed strongly regular triple with zero feedthrough, then 
		\begin{align*}
			y(t,x,v)= C_\Lambda z(t,x,v),
		\end{align*}
		for all $ x\in \mathsf{X}$, $v\in L^p_{loc}(\R_+;\mathsf{U})$, and a.e. $t\ge 0$. Furthermore $\mathbf{H}$ is given by the following familiar looking formula 
		\begin{align}\label{trans-fct}
			\mathbf{H}(\la)=\mathsf{C}_\Lambda R(\la,\mathsf{A}_{-1})\mathsf{B}, \qquad \forall\ \text{Re}\la >\omega_0(\mathsf{A}),
		\end{align}
		and the extended operator $\F$ is given by 
		\begin{align*}
			(\F v)(t)=\mathsf{C}_\Lambda \Phi_t v,  
		\end{align*}
		for all $v\in L^p_{loc}(0,\tau;\mathsf{U})$ and a.e. $t\ge 0$. In this context, $\F$ is called the extended input-output map of  $(\mathsf{A},\mathsf{B},\mathsf{C})$. For further details; see, e.g., \cite{WR}.

		The following result demonstrates a key property of $L^p$-well-posed positive systems.
		\begin{proposition}\label{Main1}
			Let $(\mathsf{A},\mathsf{B},\mathsf{C})$ be a positive $L^p$-well-posed triple with $1\le p<\infty$. Then, the following assertions are equivalent:
			\begin{itemize}
				\item[\emph{(i)}] $(\mathsf{A},\mathsf{B},\mathsf{C}) $ is a positive $L^p$-well-posed strongly regular triple with zero feedthrough.
				\item[\emph{(ii)}] $(\mathsf{A},\mathsf{B},\mathsf{C}) $ is a positive $L^p$-well-posed weakly regular triple with zero feedthrough.
			\end{itemize}
		\end{proposition}
		\begin{proof}
			The proof of this result follows from \cite[Thm. 4.1]{El}; for completeness, we provide the details here.
			
			The implication $\emph{(i)} \implies \emph{(ii)}$ is straightforward. Conversely, assume that item \emph{(ii)} holds, and let $\mathsf{C}_\Lambda$ be the Yosida extension of $\mathsf{C}$ w.r.t. $\mathsf{A}$. Then, according to \cite[Thm. 5.6.5]{Staf}, we have 
			\begin{align*}
				\langle\mathbf{H}(\la)v,y^*\rangle= \langle \mathsf{C}_\Lambda R(\la,\mathsf{A}_{-1})\mathsf{B}v,y^*\rangle, \; \forall v\in \mathsf{U},\, y^*\in \mathsf{U}'.
			\end{align*}
			On the other hand, for $s(\mathsf{A})<\lambda_1\le \la_2 $, we have 
			\begin{align*}
				R(\lambda_1,\mathsf{A})-R(\la_2,\mathsf{A})=(\la_2-\lambda_1)R(\lambda_1,\mathsf{A})R(\la_2,\mathsf{A})\ge 0,
			\end{align*}
			since $\mathsf{A}$ is resolvent positive. Multiplying this inequality by $\mathsf{C}_\Lambda$ from the left and by $\mathsf{B}$ from the right, using the positivity of $\mathsf{C}_\Lambda$ and $\mathsf{B}$, one obtains
			\begin{align}\label{deac}
				0\le \mathsf{C}_\Lambda R(\la_2,\mathsf{A}_{-1})\mathsf{B}\le  \mathsf{C}_\Lambda R(\la_1,\mathsf{A}_{-1})\mathsf{B}.
			\end{align}
			Now, let $v\in \mathsf{U}_+$ and $(\la_n)_{n\in \mathbb{N}}$ be an increasing sequence in $(s(\mathsf{A}),+\infty)$ tending to $+\infty$. It follows from \eqref{deac} and \emph{(ii)} that the sequence $(\mathsf{C}_\Lambda R(\la_n,\mathsf{A}_{-1})\mathsf{B}v)_{n\in \mathbb{N}}$ decreases monotonically and weakly converges to $0$. Thus, using Dini’s theorem in Banach lattices (see, e.g., \cite[Prop. 10.9]{BFR}),
			\begin{align*}
				\lim_{n\to +\infty } \mathsf{C}_\Lambda R(\la_n,\mathsf{A}_{-1})\mathsf{B}v=0, \qquad \forall v\in \mathsf{U}.
			\end{align*}
			Hence, the triple $(\mathsf{A},\mathsf{B},\mathsf{C}) $ is strongly regular. The proof is complete.
		\end{proof}
		
		\begin{remark}
			The above result reveals a key feature of positive $L^p$-well-posed linear systems: weak and strong regularity are equivalent. This property is especially valuable for studying control properties in linear infinite-dimensional systems, including feedback stabilization and optimal control problems.
		\end{remark}
		
		Consequently, we obtain the following results.
		\begin{corollary}
			Every positive $L^1$-well-posed triple on Banach lattices is strongly regular.
		\end{corollary}
		\begin{proof}
			The proof follows from \cite[Thm. 5.6.6]{Staf} and Proposition \ref{Main1}.
		\end{proof}
		
		\begin{corollary}
			The dual of a positive $ L^p $-well-posed strongly regular triple $(\mathsf{A},\mathsf{B},\mathsf{C})$ (with $1<p<\infty$) on reflexive Banach lattices is also a positive $ L^p $-well-posed strongly regular triple.
		\end{corollary}
		\begin{proof}
			We know that the dual of a positive $ L^p $-well-posed strongly regular triple $(\mathsf{A},\mathsf{B},\mathsf{C})$ is a positive $ L^p $-well-posed weakly regular triple. Now, Proposition \ref{Main1} yields the claim.
		\end{proof}
		
		We conclude this subsection with a generation result, which is a variant of the Weiss-Staffans perturbation theorem \cite[Thm. 6.1, Thm. 7.2 and Prop. 7.10]{WF} and \cite[Thm. 7.1.2 and Thm. 7.5.3]{Staf}. To this end, we first recall the notion of admissible positive feedback operators \cite{El}.
		\begin{definition}\label{Definition-admissible}
			Let $(\mathsf{A},\mathsf{B},\mathsf{C})$ be a positive $L^p$-well-posed triple with $1\le p<\infty$. We say that the identity $I_{\mathsf{U}}$ is an admissible positive feedback operator for $(\mathsf{A},\mathsf{B},\mathsf{C}) $ if $I_{\mathsf{U}}-\mathbb{F}$ has a positive inverse in $\calL(L^p(0,\t_0;\mathsf{U}))$ for some $\t_0>0$.
		\end{definition}
		
		We state the following characterizations of admissible positive feedback operators.
		\begin{lemma}\label{admissibil-feedback}%>\omega_0(\mathsf{A})
			Let $(\mathsf{A},\mathsf{B},\mathsf{C})$ be a positive $L^p$-well-posed triple (on $\mathsf{U},\mathsf{X},\mathsf{U}$) with $1\le p<\infty$. For $\alpha\in \R$, we define
			\begin{align*}%\label{input-output-operator}
				\F_\alpha=e_{-\alpha}\F e_{\alpha},
			\end{align*}
			where the operator $e_{\alpha}$ on $L^p_{loc}(\R_+,\mathsf{U})$ is defined by $(e_{\alpha}v)(t)=e^{\alpha t}v(t)$ for all $t\ge 0$. Then, the following conditions are equivalent:
			\begin{itemize}
				\item[(\emph{i})] $I_{\mathsf{U}}$ is an admissible positive feedback operator for $(\mathsf{A},\mathsf{B},\mathsf{C}) $.
				\item[(\emph{ii})] $r(\F)<1$.
				\item[(\emph{iii})]  $r(\F_{\alpha_0})<1$ for some $\alpha_0>\omega_0(\mathsf{A})$.
			\end{itemize}
		\end{lemma}
		\begin{proof}
			Since $\F$ is a positive operator, then $I_{\mathsf{U}}-\mathbb{F}$ has a positive inverse in $\calL(L^p(0,\t_0;\mathsf{U}))$ if and only if $r(\F)<1$. This proves (i) $\iff$ (ii). Similarly, using the positivity of $\F$ and $e_\alpha$, it follows from \cite[Thm. 7.1.8-(v)]{Staf} that $I_{\mathsf{U}}$ is an admissible positive feedback operator for $(\mathsf{A},\mathsf{B},\mathsf{C})$ if and only if $r(\F_{\alpha_0})<1$ for some $\alpha_0>\omega_0(\mathsf{A})$. This establishes (\emph{i}) $\iff$ (\emph{iii}). This completes the proof. 
		\end{proof}
		
		\begin{remark}\label{Remark-feedabck-admissibility}
			Note that when $\mathsf{X}$ and $\mathsf{U}$ are Hilbert lattices, the positivity of $\mathbf{H}$ implies (by \cite[Prop. 3.6]{WF}) that $I_{\mathsf{U}}$ is an admissible positive feedback operator for $(\mathsf{A},\mathsf{B},\mathsf{C}) $ if and only if  $r(\mathbf{H}(\la_0))<1$ for some $\la_0> \omega_0(\mathsf{A})$.
		\end{remark}	
		
		We now state and prove a version of the Weiss-Staffans perturbation theorem for positive systems.
		\begin{theorem}\label{Thoerem}
			Let $(\mathsf{A},\mathsf{B},\mathsf{C})$ be a positive $L^p$-well-posed weakly regular triple with zero feedthrough and $I_{\mathsf{U}}$ as an admissible positive feedback operator. Define 
			\begin{align*}
				\mathsf{W}_{-1}&:=(\la I_X-\tilde{\mathsf{A}}_{-1})D(\mathsf{C}_\Lambda), \qquad \la\in \rho(\mathsf{A}),
			\end{align*}
			where $\tilde{\mathsf{A}}_{-1}$ is the restriction of $\mathsf{A}_{-1}$ to $D(\mathsf{C}_\Lambda)$. Then, the operator
			\begin{align*}
				\mathsf{A}^{cl}:=\tilde{\mathsf{A}}_{-1}+\mathsf{B}\mathsf{C}_\Lambda,\qquad D(\mathsf{A}^{cl}):=\left\{x\in D(\mathsf{C}_\Lambda):\; (\tilde{\mathsf{A}}_{-1}+\mathsf{B}\mathsf{C}_\Lambda)\in \mathsf{X}\right\},
			\end{align*}
			generates a positive C$_0$-semigroup $\mathsf{T}^{cl}:=(\mathsf{T}^{cl}(t))_{t\ge 0}$ on $ \mathsf{X} $ given by 
			\begin{align*}%\label{variation}
				\mathsf{T}^{cl}(t)x=\mathsf{T}(t)x+\int_{0}^{t}\mathsf{T}_{-1}(t-s)\mathsf{B}\mathsf{C}_{\Lambda}\mathsf{T}^{cl}(s)xds,
			\end{align*}
			for all $t\ge 0$ and $x\in \mathsf{X}$.	Moreover, $(\mathsf{A}^{cl},\mathsf{B}^{cl},\mathsf{C}^{cl}_\Lambda)$ is a positive $L^p$-well-posed strongly regular triple with
			\begin{eqnarray}\label{closed-loop-operators}
				\mathsf{C}^{cl}_\Lambda=\mathsf{C}_\Lambda\in \calL(D(\mathsf{C}_\Lambda),\mathsf{U}), \qquad \mathsf{B}^{cl}=\mathsf{B}\in \calL(\mathsf{U},\mathsf{W}_{-1}).
			\end{eqnarray}
			In addition, we have $s(\mathsf{A}^{cl})\ge s(\mathsf{A})$.
		\end{theorem} 
		\begin{proof}
			By Proposition \ref{Main1}, the triple $(\mathsf{A},\mathsf{B},\mathsf{C})$ is strongly regular. Therefore, according to \cite[Thm. 7.5.3-(iii)]{Staf},  $(\mathsf{A}^{cl},\mathsf{B}^{cl},\mathsf{C}^{cl}_\Lambda)$ is an $L^p$-well-posed strongly regular triple. To complete the proof, we only need to verify that $(\mathsf{A}^{cl},\mathsf{B}^{cl},\mathsf{C}^{cl}_\Lambda)$ is positive. In fact, let $\la\in \rho(\mathsf{A}) \cap \rho(\mathsf{A}^{cl})$ (since $\rho(\mathsf{A})\cap \rho(\mathsf{A}^{cl})\neq \emptyset$). By \cite[Thm. 7.4.7-(viii)]{Staf}, the operator $I_{\mathsf{U}}-\mathbf{H}(\la_0)$ is invertible, and 
			\begin{eqnarray*}
				R(\la_0,\mathsf{A}^{cl})=\left(I_{\mathsf{X}}+R(\la_0,\tilde{\mathsf{A}}_{-1})\mathsf{B}(I_{\mathsf{U}}-\mathbf{H}(\la_0))^{-1}\mathsf{C}\right)R(\la_0,\mathsf{A}).
			\end{eqnarray*}
			Using the definition of the spectral radius, the monotonicity of the operator norm, and the positivity and monotonic decrease of the family $(\mathbf{H}(\la))_{\la >\omega_0(\mathsf{A})}$, we obtain  
			\begin{align*}
				r(\mathbf{H}(\la))\le r(\mathbf{H}(\la_0))<1,\qquad \forall \la \ge \la_0>\omega_0(\mathsf{A}).
			\end{align*}
			Thus, for all $\la\ge \la_0$,
			\begin{align}\label{resolvent-closed-loop}
				R(\la,\mathsf{A}^{cl})=R(\la,\mathsf{A})+\sum_{k\ge 0}R(\la,\tilde{\mathsf{A}}_{-1})\mathsf{B} \mathbf{H}^k(\la)\mathsf{C} R(\la,\mathsf{A})
				\ge R(\la,\mathsf{A})\ge  0,
			\end{align}
			which implies that $\mathsf{A}^{cl}$ is a resolvent positive operator. Hence, the semigroup $\mathsf{T}^{cl}$ is positive. Furthermore, form \eqref{closed-loop-operators}, it follows that the operators $\mathsf{B}^{cl}$ and $\mathsf{C}_\Lambda^{cl}$ are positive. Thus, $(\mathsf{A}^{cl},\mathsf{B}^{cl},\mathsf{C}_\Lambda^{cl})$ is a positive triple.
			
			Now, assume that $s(\mathsf{A}^{cl}),s(\mathsf{A})>-\infty$ (otherwise, there is nothing to prove). On the one hand, by \eqref{resolvent-closed-loop}, we have
			\begin{align*}
				\Vert R(\la,\mathsf{A}^{cl})\Vert \ge \Vert R(\la,\mathsf{A})\Vert,
			\end{align*}
			for all $\la>\max\{s(\mathsf{A}),s(\mathsf{A}^{cl})\}$. On the other hand, by \cite[Cor. 12.9]{BFR}, $s(\mathsf{A})\in \sigma(\mathsf{A})$, and thus
			\begin{align*}
				\lim_{ \la \to s(\mathsf{A})} \Vert R(\la,\mathsf{A})\Vert=\infty. 
			\end{align*}
			Therefore, $s(\mathsf{A}^{cl})\ge s(\mathsf{A})$. This completes the proof. 
		\end{proof}
		
		\begin{remark}
			Theorem \ref{Thoerem} provides mild sufficient conditions under which the positivity of the triple $(\mathsf{A},\mathsf{B},\mathsf{C})$ is preserved under static feedback represented by a linear bounded positive operator. For instance, Theorem \ref{Thoerem} can significantly simplify the positivity analysis of the semigroup conducted in \cite[Sect. 4]{BFR1}. More precisely, assumption (4.2) in \cite[Thm. 4.2]{BFR1} can be dropped. In fact, as shown in the proof of Theorem \ref{Thoerem}, the $L^p$-well-posedness of both the open-loop and closed-loop systems ensures that $\rho(\mathsf{A}) \cap \rho(\mathsf{A}^{cl}) \neq \emptyset$, cf. \cite[Thm. 7.4.7]{Staf}. Moreover, using the positivity of the transfer function $\mathbf{H}$, it follows from \cite[Thm. 7.4.7-(viii)]{Staf} that there exists $\lambda_0 > \omega_0(\mathsf{A})$ such that $r(\mathbf{H}(\lambda_0)) < 1$.
		\end{remark}

		\section{Main results}\label{Sec:3}
		In this section, we present and prove the main results of the paper. We start with the well-posedness of \eqref{Eq1}.
		\subsection{Well-posedness}\label{Sec:3.1}
		To state our results rigorously, we need to introduce some notations and concepts. We consider the operator $Q:D(Q)\to Y_p$ defined by 
		\begin{align*}
			Q\varphi:=\partial_\theta \varphi, \qquad \varphi\in D(Q):=\ker(\delta_0),
		\end{align*}
		where $\delta_0:W^{1,p}(-r,0;V)\to V$ is the Dirac mass at $0$, defined by $\delta_0\varphi:=\varphi(0)$. It is well-known (see, e.g., \cite{BFR}) that the operator $Q$ generates the nilpotent left shift semigroup $\S:=(S(t))_{t\ge 0}$ on $Y_p$ given by 
		\begin{align*}%\label{S-sg}
			(S(t)\varphi)(\theta)=\begin{cases} 0,& -t\le \theta\le 0,\cr \varphi(t+\theta),& -r\le \theta\le -t,\end{cases}
		\end{align*}
		for all $t\ge 0$, $\varphi\in Y_p$, and a.e. $\theta\in [-r,0]$. We also introduce the operator $\varepsilon_\la:V\to Y_p$ defined by $(\varepsilon_\la v)(\theta):=e^{\la \theta}v$ for $v\in V$, $\la \in \C$, and $\theta\in [-r,0]$.
		\vspace{.1cm}
		
		We state the following result.
		\begin{lemma}\label{regular-shift-triple}
			Let Assumption \ref{Assp22} be satisfied and the delay operator $L$, defined by \eqref{dealy-operator}, be positive. We select the operators
			\begin{eqnarray}\label{beta}
				\Theta:= (-Q_{-1})\varepsilon_0\in \calL(V, (Y_p)_{-1}^Q) \qquad \L:=L\vert_{D(Q)}.
			\end{eqnarray} 
			Then, $(Q,\Theta,\L)$ is a positive $L^p$-well-posed uniformly regular triple with zero feedthrough and input-maps $\Phi^{Q,\Theta}_t$ given by
			\begin{align}\label{control-maps-beta}
				(\Phi^{Q,\Theta}_t \psi)(\theta)=\begin{cases} \psi(t+\theta),& -t\le \theta \le 0,\cr 0,& -r\le \theta\le -t,\end{cases}
			\end{align}
			for all $t\ge 0$, $\psi\in L^p(\R^+,V)$, and a.e. $\theta\in [-r,0]$.
		\end{lemma}
		\begin{proof} 
			The proof of this lemma can be deduced from the abstract result of \cite[Thm. 3]{HIR}. However, since we are particularly interested in the positivity and uniform regularity of the triple $(Q,\Theta,\L)$, we give here the proof. Indeed, using H\"{o}lder's inequality and Fubini's theorem, one obtains that
			\begin{align*}
				\int_{0}^{\alpha} \Vert \mathbb{L} S(t)\varphi\Vert^p dt \le   {\rm Var}(\eta;-r,0)^p\Vert\varphi\Vert_{Y_p}^p,
			\end{align*}
			for any $0<\alpha<r$ and any $ \varphi\in D_+(Q)$. Therefore, $(\mathbb{L},Q)$ is positive $L^p$-admissible. In addition, using the injectivity of the Laplace transform, we obtain that the input-map $\Phi^{Q,\Theta}_\cdot$ is given by \eqref{control-maps-beta}. The positivity of $\Phi^{Q,\Theta}_t$, for every $t\ge 0$, further confirms that $(Q,\Theta)$ is positive $L^p$-admissible. 
			
			Now, for $0<\alpha<r$, we define the input-output map associated with $(Q,\Theta,\L)$:
			\begin{align}\label{Shift-estimate}
				(\mathbb{F}^{Q,\Theta,\L}\psi)(t):=L\Phi^{Q,\Theta}_t\psi,\quad t\in[0,\alpha], % $\alpha>r>-\theta$ implies $\alpha+\theta>0$.
			\end{align}
			for any $\psi\in W^{1,p}_{0,+}(0,\alpha;V)$. By the H\"{o}lder's inequality and Fubini’s theorem, we obtain
			\begin{align*}
				\int_{0}^{\alpha} \Vert (\mathbb{F}^{Q,\Theta,\L}\psi)(t)\Vert^p dt\le  {\rm Var}(\eta;-\alpha,0)^p \Vert \psi\Vert_{L^p(0,\alpha; V)}^p,
			\end{align*}
			for any $\psi\in W^{1,p}_{0,+}(0,\alpha;V)$. Thus, the triple $(Q,\Theta,\L)$ satisfies the conditions (i)-(iii) in Proposition \ref{ABC-sigma} and therefore $(Q,\Theta,\L)$ is a positive $L^p$-well-posed triple on $V,Y_p,V$. Furthermore, one has
			\begin{align*}
				\left\Vert\mathbb{F}^{Q,\Theta,\L}(\1_{\mathbb{R}_{+}} \cdot v)\right\Vert_{L^p(0,\alpha; V)}  
				&\le  \alpha^{\frac{1}{p}} {\rm Var}(\eta;-\alpha,0) \Vert v\Vert_{V},
			\end{align*}
			for all $v\in V$, where ${\rm Var}(\eta;-\alpha,0)\to 0$ as $\alpha\downarrow 0$ (cf. Remark \ref{further1}). By \cite[(ii)-Thm 5.6.7]{Staf}, we conclude that $(Q,\Theta,\L)$ is a positive $L^p$-well-posed uniformly regular triple with zero feedthrough.  The proof is complete.
		\end{proof}
		
		We can now state and prove the first main result of this paper.
		\begin{theorem}\label{Main1-23}
			Let Assumptions \ref{Assp22} and \ref{Assp11} be satisfied. Suppose that $\T$ is positive, $P$ is positive, and $\D_\la $ is positive for all sufficiently large $\la\in \R$.  Assume the following:
			\begin{enumerate}
				\item[(a)] $L$ is positive;
				\item[(b)] $ (A,\B,\P) $ is a positive $L^p$-well-posed weakly regular triple with zero feedthrough, where $\P:=P\vert_{D(A)}$;
				\item[(c)] $r(\mathbb{F}^{Q,\Theta,\L}\F^{A,\B,\P})<1$, where $\mathbb{F}^{Q,\Theta,\L}$ and $\F^{A,\B,\P}$ are the input-output maps of $(Q,\Theta,\L)$ and $ (A,\B,\P) $, respectively.
			\end{enumerate}
			Then, the operator $\mathfrak{A}$, defined by \eqref{cauchy-operator}, generates a positive C$_0$-semigroup $(\mathfrak{T}(t))_{t\ge 0}$ on $\mathscr{X}_p$. Moreover, the boundary delay equation \eqref{Eq1}  has a unique positive mild solution $z(\cdot)\in \calC(\R_+,X)$ given by
			\begin{align*}
				z(t,f,\varphi)=\left[\mathfrak{T}(t)(f,\varphi)^\top\right]_1,
			\end{align*}
			for all $t\ge 0$ and $(f,\varphi)^\top\in \mathscr{X}_p$, where $\left[\mathfrak{T}(t)(f,\varphi)^\top\right]_1$ denotes the first component of $\mathfrak{T}(t)(f,\varphi)^\top$.
		\end{theorem} 
		\begin{proof}%[Proof of Theorem \ref{Main1-23}]
			First, we consider the operator $\mathscr{A}_{m}:D(\mathscr{A}_{m})\to \mathscr{X}_p$ defined by
			\begin{eqnarray*}%\label{A_M}
				\mathscr{A}_{m}:=\begin{pmatrix}
					A_m & 0 \\
					0 & \partial_\theta \\
				\end{pmatrix}, \quad D(\mathscr{A}_{m}):=D(A_m)\times W^{1,p}(-r,0;V).
			\end{eqnarray*}
			In addition, we select the operators matrices $\mathscr{G},\mathscr{P}:D(\mathscr{A}_m)\to V\times V:= \mathscr{V}$ defined by
			\begin{align*}%\label{G,Ga}
				\begin{array}{lll}
					\mathscr{G}=\begin{pmatrix}
						G & 0 \\
						0 & \delta_{0} \\
					\end{pmatrix},\qquad
					\mathscr{P}= \begin{pmatrix}
						0 & 	L \\
						P & 0 \\
					\end{pmatrix}.
				\end{array}
			\end{align*}
			With these notations, the domain of $\mathfrak{A}$ can be rewritten as
			\begin{eqnarray*}%\label{Big-A}
				D(\mathfrak{A})=\big\{(f,\varphi)^\top\in D(\mathscr{A}_m):\: \mathscr{G}(f,\varphi)^\top=\mathscr{P}(f,\varphi)^\top \big\}.
			\end{eqnarray*}
			To prove our claim, we shall use Theorem \ref{theorem-ABC}. In fact, we set
			\begin{eqnarray*}
				\mathscr{A}:=\mathscr{A}_m,\qquad D(\mathscr{A}):=\ker(\mathscr{G})=\ker(G)\times \ker(\delta_0).
			\end{eqnarray*}
			By  Assumption \ref{Assp11}-(i), the operator 
			\begin{align*}
				\mathscr{A}={\rm diag}(A,Q),\qquad D(\mathscr{A})=\ker(G)\times \ker(\delta_0),
			\end{align*}
			generates a positive C$_0$-semigroup $\mathscr{T}:=(\mathscr{T}(t))_{t\ge 0}$ on $\mathscr{X}_p$, given by $\mathscr{T}:={\rm diag}(\T,\S)$. Moreover, Assumption \ref{Assp11}-(ii) implies that the operator $\mathscr{G}$ is onto. Then, one can readily show that the Dirichlet operator associated with $(\mathscr{A}_m,\mathscr{G})$ is given by
			\begin{align*}
				\mathscr{D}_\la={\rm diag}(\D_\la,\varepsilon_\la)\in \calL\left(\mathscr{V},\mathscr{X}_p\right), \qquad \la >s(A),
			\end{align*}
			where $\D_\la$ and $\varepsilon_\la$ are the Dirichlet operators associated with $(A_m,G)$ and $(\partial_\theta,\delta_{0})$, respectively. We then define 
			\begin{align*}
				\mathscr{B}:=(\la I_{\mathscr{X}_p}-\mathscr{A}_{-1})\mathscr{D}_\la\in \calL\left(\mathscr{V},(\mathscr{X}_p)_{-1}^\mathscr{A}\right), \quad \la> s(A).
			\end{align*}
			Thus, $\mathscr{B}={\rm diag}(\B,\Theta)$. Since $(A,\B)$ and $(Q,\Theta)$ are both positive $L^p$-admissible, then $(\mathscr{A},\mathscr{B})$ is also positive $L^p$-admissible. On the other hand, we have
			\begin{eqnarray*}
				\mathscr{C}:=\mathscr{P}\vert_{D(\mathscr{A})}=\begin{pmatrix}
					0 & \L \\
					\P & 0
				\end{pmatrix}.
			\end{eqnarray*}
			Using assumption (b) and according to Lemma \ref{regular-shift-triple}, we obtain that $(\mathscr{C},\mathscr{A})$ is positive $L^p$-admissible. Its Yosida extension is given by
			\begin{eqnarray*}%\label{C-extension}
				\mathscr{C}_\Lambda=\begin{pmatrix}
					0 & \L_{\Lambda}\\
					\P_{\Lambda}& 0
				\end{pmatrix}, \qquad D(\mathscr{C}_\Lambda)=D(\P_{\Lambda})\times D(\L_{\Lambda}),
			\end{eqnarray*}
			where $\P_{\Lambda}$ and $\L_{\Lambda}$ are the Yosida extensions of $\P$ and $\L$ with respect to $A$ and $Q$, respectively. 
			\vspace{.1cm}
			
			Now, for $\alpha>r$, we define the input-output map associated with $(\mathscr{A},\mathscr{B},\mathscr{C})$:
			\begin{eqnarray*}%\label{Shift-estimate}
				(\mathbb{F}(u,v)^\top)(t):=\mathscr{P}\Phi_t(u,v)^\top=\begin{pmatrix}
					L\Phi_{t}^{Q,\Theta}v\\
					P\Phi_t^{A,\B}u
				\end{pmatrix},\quad t\in[0,\alpha],
			\end{eqnarray*} 
			for all $(u,v)^\top\in W^{1,p}_{0,+}(0,\alpha; \mathscr{V})$, where $\Phi_{t}^{Q,\Theta}$ and $\Phi_{t}^{A,\B}$ are the input-maps of $(Q,\Theta)$ and $(A,\B)$, respectively. Using the assumption (b) and Lemma \ref{regular-shift-triple}, it follows that $\F\in \calL( L^p_+(0,\alpha; \mathscr{V}))$ and 
			\begin{align}\label{Big-F}
				\F =\begin{pmatrix}
					0& \F^{Q,\Theta,\L}\\
					\F^{A,\B,\P}& 0
				\end{pmatrix},
			\end{align} 
			where $\F^{Q,\Theta,\L}$ and $\F^{A,\B,\P}$ are the input-output maps of the triples $(Q,\Theta,\L)$ and $(A,\B,\P)$, respectively. Thus, $(\mathscr{A},\mathscr{B},\mathscr{C})$ is a positive $L^p$-well-posed triple on $\mathscr{V},\mathscr{X}_p,\mathscr{V}$. Taking the Laplace transform in \eqref{Big-F}, we obtain 
			\begin{eqnarray*}
				\mathbf{H}(\la)=\begin{pmatrix}
					0&\mathbf{H}^{Q,\Theta,\L}(\la)\\
					\mathbf{H}^{A,\B,\P}(\la)&0
				\end{pmatrix},\qquad \forall\la>\max\{\omega_0(A),\omega_0(Q)\},
			\end{eqnarray*}
			where $\mathbf{H}^{Q,\Theta,\L}$ and $\mathbf{H}^{A,\B,\P}$ are the transfer functions of the triples $(Q,\Theta,\L)$ and $(A,\B,\P)$, respectively. Therefore, 
			\begin{eqnarray*}
				\lim_{\R\ni \la\to +\infty }\langle\mathbf{H}(\la)g,h^*\rangle=0,\qquad \forall g\in \mathscr{V},\,h^*\in \mathscr{V}',
			\end{eqnarray*}
			since $(Q,\Theta,\L)$ and $(A,\B,\P)$ are weakly regular triples with zero feedthrough. We conclude that $(\mathscr{A},\mathscr{B},\mathscr{C})$ is a positive $L^p$-well-posed weakly regular triple with zero feedthrough.
			
			Next, we prove that the identity $I_{\mathscr{V}}$ is an admissible positive feedback operator for $(\mathscr{A},\mathscr{B},\mathscr{C})$. In fact, we have 
			\begin{align*}
				I_{\mathscr{V}}-\F&= \begin{pmatrix}
					I_V& -\F^{Q,\Theta,\L}\\
					-\F^{A,\B,\P}& I_V
				\end{pmatrix}\\&=\begin{pmatrix}
					I_{V}-\F^{Q,\Theta,\L}\F^{A,\B,\P}&	-\F^{Q,\Theta,\L}\\
					0 & I_{V}
				\end{pmatrix}\begin{pmatrix}
					I_{V} &	0\\
					-\F^{A,\B,\P} & I_{V} 
				\end{pmatrix}.
			\end{align*} 
			Thus, according to Lemma \ref{admissibil-feedback}, $I_{\mathscr{V}}$ is an admissible positive feedback operator if and only if $r(\F^{Q,\Theta,\L}\F^{A,\B,\P})<1$, which holds by assumption (c) . Hence, by Theorem \ref{theorem-ABC}, $\mathfrak{A}$ generates a positive C$_0$-smeigroup $\mathfrak{T}$ on $\mathscr{X}_p$ given by 
			%{\footnotesize
				\begin{eqnarray*}
					\mathfrak{T}(t)(f,\varphi)^\top=\mathscr{T}(t)(f,\varphi)^\top+\int_{0}^{t}\mathscr{T}_{-1}(t-s) \mathscr{B}\mathscr{C}_\Lambda\mathfrak{T}(s)(f,\varphi)^\top ds,
				\end{eqnarray*}
				for all $t\ge 0$ and $(f,\varphi)^\top\in \mathscr{X}_p$. Moreover, the abstract boundary delay equation \eqref{Eq1} has a unique positive mild solution 
				$z(\cdot):[0,+\infty)\to X$ given by
				\begin{align*}
					z(t,f,\varphi)=\left[\mathfrak{T}(t)(f,\varphi)^\top\right]_1,
				\end{align*}
				for all $t\ge 0$ and $(f,\varphi)^\top\in \mathscr{X}_p$. Thus, 
				%	{\footnotesize
					\begin{eqnarray}\label{variation1}
						%			\begin{array}{lll}
							z(t,f,\varphi)=T(t)f+ \int_{0}^{t}T_{-1}(t-s)\B\L_{\Lambda}S(s)\varphi ds+ \int_{0}^{t}T_{-1}(t-s)\B\L_{\Lambda}\int_{0}^{s}S_{-1}(s-\sigma)\Theta \P_\Lambda z(\sigma,f,\varphi)d\sigma ds,%\nonumber
							%		\end{array}
					\end{eqnarray}
					for all $t\geq 0$ and $(f,\varphi)^\top\in \mathscr{X}_p$, where $\L_{\Lambda}$ and $\P_\Lambda$ are the Yosida extensions of $\mathbb{L}$ and $\mathbb{P}$ w.r.t $Q$ and $A$, respectively. This completes the proof.
				\end{proof}

					\begin{corollary}\label{Cor.11}
						Let the assumptions of Theorem \ref{Main1-23} be satisfied, whereby (c) is replaced by the following:
						\begin{itemize}
							\item[(c')] For every $\alpha\ge 0$, there exists $\nu_0>0$ independent of $\alpha$ such that $\Vert \F^{A,\B,\P}\Vert_{\calL(L^p(0,\alpha;V))}\le \nu_0$.
						\end{itemize}
						Then, the abstract delay equation \eqref{Eq1}  has a unique positive mild solution $z(\cdot)\in \calC(\R_+,X)$ for all $t\ge 0$, $f\in X$, and $\varphi\in Y_p$ given by \eqref{variation1}.
					\end{corollary}
					\begin{proof}
						According to Theorem \ref{Main1-23}, it suffices to show that $r(\mathbb{F}^{Q,\Theta,\L}\F^{A,\B,\P})<1$. Indeed, using (c') and Lemma \ref{regular-shift-triple}, we obtain   
						\begin{align*}
							\int_{0}^{\alpha} \Vert (\F^{Q,\Theta,\L}\F^{A,\B,\P}g)(t)\Vert^p dt&\le	 {\rm Var}(\eta;-\alpha,0)^p \Vert \F^{A,\B,\P}g\Vert_{L^p(0,\alpha;V)}^p\\
							&	\le {\rm Var}(\eta;-\alpha,0)^p\nu_0^p\Vert g\Vert_{L^p(0,\alpha;V)}^p,
						\end{align*}
						for all $g\in L^p_+(0,\alpha;V)$ and $\alpha\in (0,r)$. Since ${\rm Var}(\eta;-\alpha,0)\to 0$ as $\alpha\downarrow 0$, then there exists $\alpha_0>0$ small enough such that
						\begin{align*}
							\Vert\F^{Q,\Theta,\L}\F^{A,\B,\P}\Vert_{\mathcal{L}(L^p_+(0,\alpha_0;V))}<1,
						\end{align*}
						as $\nu_0$ is independent of $\alpha$.
					\end{proof}
					
					\begin{remark}
						Almost all results in this section remain valid, with minor and straightforward modifications, when Banach lattices are replaced by Banach spaces. The exceptions are the following: the positivity assumptions on $\T,P,\D_\la$ and $L$ can be dropped. The weak regularity assumption on $ (A,\B,\P) $ must be strengthened to strong regularity. This is because the converse part of Proposition \ref{Main1} is a consequence of Dini theorem, which has no counterpart for Banach spaces. The assumption (c) in Theorem \ref{Main1-23} has to be replaced by $1\in \rho(\mathbb{F}^{Q,\Theta,\L}\F^{A,\B,\P})$. With these modifications, we apply \cite[Thm. 4.1]{HMR} in the proof of Theorem \ref{Main1-23}.
					\end{remark}
					
					\subsection{Stability}\label{Sec:3.2}
					In this subsection, we derive necessary and sufficient conditions for the exponential stability of \eqref{Eq1}. Building on the reformulation presented in Sections \ref{Sec:1.3} and \ref{Sec:3.1}, and according to Definition \ref{Def.A}, the abstract delay boundary problem \eqref{Eq1} is \emph{strongly exponentially stable} if there exists $\alpha>0$ such that
					\begin{align*}
						\left\Vert \mathfrak{T}(t)(f,\varphi)^\top\right\Vert_{\mathscr{X}_p}\le \Upsilon_{f,\varphi} \,e^{-\alpha t}\left\Vert (f,\varphi)^\top\right\Vert_{ \mathscr{X}_p},
					\end{align*}
					for all $(f,\varphi)^\top\in D(\mathfrak{A})$, $t\ge 0$, and a constant $\Upsilon_{f,\varphi}>0$. 
					\vspace{.1cm}
					
					Our second main result, stated in the following theorem, provides a spectral condition for the exponential stability of \eqref{Eq1}.
					\begin{theorem}\label{Main1-24}
						Let the assumptions of Theorem \ref{Main1-23} be satisfied. Define
						\begin{align*}
							\rho(A)\ni\la \mapsto \mathbb{H}(\la) := P\D_\la L(\varepsilon_\la)\in \calL(V).
						\end{align*}
						Then, $\mathfrak{T}$ is strongly exponential stable if and only if 
						\begin{align}\label{spectral-condition}
							s(A)<0\quad {\rm and}\quad r(\mathbb{H}(0))<1.
						\end{align}
						If the above condition holds true, then there exists $\alpha>0$ such that
						\begin{align*}
							\left\Vert (z(t),x_t)^\top\right\Vert_{\mathscr{X}_p}\le \Upsilon_{f,\varphi} \,e^{-\alpha t}\left\Vert (f,\varphi)^\top\right\Vert_{ \mathscr{X}_p},
						\end{align*}
						for all $(f,\varphi)^\top\in D(\mathfrak{A})$, $t\ge 0$, and a constant $\Upsilon_{f,\varphi}>0$. 
					\end{theorem}

					\begin{proof}%[Proof of Theorem \ref{Main1-24}]
						For $\la \in \rho(\mathscr{A})$, we have 
						\begin{eqnarray*}%\label{cara2}
							I_{\mathscr{V}}-\mathscr{P} \mathscr{D}_\la=\begin{pmatrix}
								I_{V} & 0\\
								-P \D_\lambda & I_{V}-\mathbb{H}(\la)
							\end{pmatrix}\begin{pmatrix}
								I_{V} &- L(\varepsilon_\la)\\
								0& I_{V}
							\end{pmatrix}.
						\end{eqnarray*}
						Then, for such $\la$, one obtains that $r(\mathscr{P} \mathscr{D}_\la)<1$ if and only if $r(\mathbb{H}(\la))<1$. Thus, according to Theorem \ref{stability-criterion}, we obtain 
						\begin{align*}
							s(\mathfrak{A})<0&\iff s(\mathscr{A})<0\quad {\rm and }\quad r(\mathscr{P} \mathscr{D}_0)<1\\
							&\iff s(A)<0\quad {\rm and }\quad r(\mathbb{H}(0))<1,
						\end{align*}
						where we used the fact that $s(Q)=-\infty$ to obtain the second equivalence. Thus, $\mathfrak{T}$ is strongly exponential stable if and only if the spectral condition \eqref{spectral-condition} holds. This ends the proof.
					\end{proof}
					
					Using the spectral properties of positive semigroups, this result can be improved for $L^p$-spaces.
					\begin{corollary}\label{Cor1}
						Consider the setting of Theorem \ref{Main1-24}, let $X:=L^p(\Omega_1,\Sigma_1, \mu_1)$ and $V:=L^p(\Omega_2,\Sigma_2, \mu_2)$, where $p\in [1,\infty)$ and $(\Omega_i,\Sigma_i, \mu_i)$ are $\sigma$-finite measure spaces for $i=1,2$. If the spectral condition \eqref{spectral-condition} holds, then there exist $\Upsilon, \alpha>0$ such that 
						\begin{align*}
							\left\Vert (z(t),x_t)^\top\right\Vert_{\mathscr{X}_p}\le \Upsilon e^{-\alpha t}\left\Vert (f,\varphi)^\top\right\Vert_{ \mathscr{X}_p},
						\end{align*}
						for all $t\ge 0$ and $(f,\varphi)^\top\in \mathscr{X}_p$, where $\Upsilon>0$ is independent of $f$ and $\varphi$.
					\end{corollary}
					\begin{proof}
						The proof follows from Theorem \ref{Main1-24}, since in this case we have $s(\mathfrak{A})=\omega_0(\mathfrak{A})$ (see \cite{Weis}).
					\end{proof}
					
					In the following, we establish necessary and sufficient conditions for the exponential stability of two fundamental classes of boundary delay equations.
					\begin{corollary}\label{Cor2}
						Let Assumptions \ref{Assp11} and \ref{Assp22} be satisfied, and suppose conditions (b)-(c) of Theorem \ref{Main1-23} hold.  Let $\eta$ be an operator-valued increasing function, i.e., 
						\begin{align*}
							0=\eta(-r)\le \eta(\theta_1) \le \eta(\theta_2),\quad -r\le \theta_1\le \theta_2\le 0.
						\end{align*}
						Then, the boundary delay equation \eqref{Eq1} is strongly exponentially stable if and only if
						\begin{align}\label{C1}
							s(A)<0\quad {\rm and}\quad r(P\D_0\eta(0))<1.
						\end{align}
					\end{corollary}
					\begin{proof}
						Since $\eta$ is an operator-valued increasing function, it follows from the closedness of $V_+$ that the delay operator $L$ is positive. Thus, according to Theorem \ref{Main1-23}, system \eqref{Eq1} is well-posed. Furthermore, since $\eta(-r)=0$, we have $L(\varepsilon_0)=\eta(0)$, which implies that $\mathbb{H}(0)=P\D_0\eta(0)$. Therefore, by Theorem \ref{Main1-24}, the boundary delay equation \eqref{Eq1} is strongly exponentially stable if and only if condition \eqref{C1} holds. This completes the proof. 
					\end{proof}
					
					\begin{corollary}\label{Cor3}
						Let Assumptions \ref{Assp11} and \ref{Assp22} be satisfied. Assume that the assumptions (b)-(c) in Theorem \ref{Main1-23} hold. If we take $\Lambda\in \calL_+(V)$ and consider
						\begin{align*}
							\eta(\theta)=\begin{cases}
								\Lambda,& {\rm if}\;\theta=-r,\\
								0,& {\rm if}\; \theta\in (-r,0],
							\end{cases}
						\end{align*}   
						then the boundary delay equation \eqref{Eq1} is strongly exponentially stable if and only if
						\begin{align}\label{C2}
							s(A)<0\quad {\rm and}\quad r(P\D_0\Lambda)<1.
						\end{align}
					\end{corollary}
					\begin{proof}
						Notice that in this case, we have ${\rm Var}(\eta;-r,0)\le \Vert\Lambda\Vert_{\calL(V)}$, and therefore $\eta$ satisfies Assumption \ref{Assp22}. Moreover, we have $L(\varphi)=\varphi(-r)$, which implies that the delay operator $L$ is positive. Thus, according to Theorem \ref{Main1-23}, system \eqref{Eq1} is well-posed. Furthermore, by Theorem \ref{Main1-24}, the boundary delay equation \eqref{Eq1} is strongly exponentially stable if and only if condition \eqref{C2} holds. This completes the proof.
					\end{proof}
					
					\begin{remark}\label{unifrom-stability}
						By Corollary \ref{Cor1}, the results of Corollaries \ref{Cor2} and \ref{Cor3} characterize the uniform exponential stability of system \eqref{Eq1} when $X$ and $V$ are $L^p$-spaces ($1\le p<\infty$).
					\end{remark}

					\section{Application} \label{Sec:4}
			In this section, we are concerned with the exponential stability of a one-dimensional hyperbolic system of the form
			\begin{eqnarray}\label{Eq2}
				\partial_t z(t,x)+\mathbf{V}(x)\partial_xz(t,x)=\mathbf{K}(x)z(t,x),\; t\geq 0,x\in(0,\ell),
			\end{eqnarray}
			subject to the following time-delayed proportional boundary feedback law
			\begin{eqnarray}\label{Eq3}
				z(t,0)=\int_{-r}^{0}d\eta(\theta) z(t+\theta,\ell), \qquad t\geq 0,
			\end{eqnarray}
			and initial conditions
			\begin{eqnarray*}
				z(0,x)= f(x),\quad z(\theta,\ell)= \varphi(\theta), \quad x\in(0,\ell), \theta\in [-r,0].
			\end{eqnarray*}
			Here,
			\begin{itemize}
				\item[(a)] $z:\R_+\times [0,\ell] \to \R^n$ ($\ell>0$ and $n\in \N^*$) describes the system dynamics;
				\item[(b)] The matrices $\mathbf{V}:={\rm diag }(v_{i}(\cdot))_{1\le i\le n}$ and $\mathbf{K}:={\rm diag}(k_i(\cdot))_{1\le i\le n}$ are real, diagonal, and bounded, i.e., $v_{i}(\cdot),k_{i}(\cdot)\in L^{\infty}(0,\ell)$  for all $i\in\{1,\ldots,n\}$;
				\item[(c)] There exists $\nu>0$ such that $v_{i}(x)\ge \nu$ for all $i\in\{1,\ldots,n\}$ and a.e. $x\in (0,\ell)$;
				\item[(d)] The matrix function $\eta:[-r,0]\to \R^{n\times n}$ is real-valued and satisfies Assumption \ref{Assp22};
				\item[(e)] The Riemann-Stieltjes integral $\int_{-r}^{0}d\eta(\theta)$ define a positive operator.
			\end{itemize}
			
			In what follows, we denote system \eqref{Eq2} with the boundary conditions \eqref{Eq3} by $(\Sigma_{\mathsf{Hyp}})$. To investigate the well-posedness and stability of $(\Sigma_{\mathsf{Hyp}})$, we first recast the system within an abstract operator framework on suitable function spaces. To this end, we consider the following Banach space
			\begin{align*}
				X_p:= L^p(0,\ell;\R^n), \qquad 1\le p< \infty.
			\end{align*}
			On this space we introduce the operator
			\begin{eqnarray}\label{A_m}
				A_mf:= - \mathbf{V}(\cdot)\partial_x f +\mathbf{K}(\cdot)f, \;  D(A_m):= W^{1,p}(0,\ell;\R^n).
			\end{eqnarray}
			Moreover, we consider the boundary operators $G,P:D(A_m)\to \R^n$ defined as 
			\begin{eqnarray}\label{trace-operators}
				Gf=f(0), \quad Pf=f(\ell), \qquad \forall f\in D(A_m).
			\end{eqnarray}
			We now consider the Banach space
			\begin{align*}
				Y_{p}:=L^p(-r,0;\R^n),
			\end{align*}
			and select the delay operator 
			\begin{align*}
				L(\varphi):=\int_{-r}^{0} d\eta(\theta) \varphi(\theta), \qquad \varphi\in  W^{1,p}(-r,0;\R^n).
			\end{align*}
			
			With these notations, the hyperbolic system $(\Sigma_{\mathsf{Hyp}})$ can be written in the preliminary abstract form
			\begin{align*}%\label{Eq2}
				\begin{cases}
					\partial_t z(t)=A_mz(t),&t> 0,\cr  
					Gz(t)= L(\tilde{z}_t), & t> 0,\cr
					\tilde{z}(t)=Pz(t),& t>0,\\
					z(0)=f,\; \tilde{z}_0=\varphi, & 
				\end{cases}
			\end{align*}%that 
			where $\tilde{z}_t$ is the history segment of the state value $z(t)$ at the endpoint $\ell$. Before proceeding, we make the following remark.
			\begin{remark}
				We note that \cite{Zhang} showed that for the (uniform) delay-free case (i.e., $\mathbf{K}\equiv 0$ and $\mathbf{V},\eta$ are positive constant matrices), $(\Sigma_{\mathsf{Hyp}})$ is uniformly exponentially stable if and only if $r(\eta)<1$, provided the condition that the inverse velocities $v_i$ are commensurate. \cite{BBH} extended this, proving that the same stability criterion holds without the commensurate velocity assumptions. Their work also examines a class of linear hyperbolic systems with time-delayed boundary conditions bounded with respect to the state variable; see \cite[Syst. 5.1]{BBH}. As a consequence, the proposed boundary coupling can be impractical in real-world applications. For example, in freeway traffic control, this would require impractical continuous sensing and actuation throughout the freeway network.  In contrast, proportional boundary feedback control offers a more feasible alternative using naturally occurring control points such as on-ramp metering and variable message signs \cite{Zhang}.
			\end{remark}
			%\vspace{.1cm}
			
			Under the hypothesis (b) and (c), the following functions are well-defined for all $x\in [0,\ell]$ and $i\in\{1,\ldots,n\}$:
			\begin{eqnarray*}%\label{S3.functions}
				\tau_{i}(x):=\int_{0}^{x}\dfrac{dx}{v_{i}(x)},\quad
				\xi_{i}(x):=\int_{0}^{x} \dfrac{k_{i}(x)}{v_{i}(x)}dx.
			\end{eqnarray*}
			Moreover, the hypothesis (c) ensures that each $\tau_i(\cdot)$ is strictly increasing and continuous, and therefore invertible. Let $\t_i^{-1}(\cdot)$ denote its inverse function. For any initial position $x_0\in [0,\ell]$, we define the characteristic trajectory
			\begin{align*}%\label{s-tilde} 
				\pi_i(t)=\t_i^{-1}(\t_i(x_0)-t),\qquad t\in [0,\t_i(\ell)], 
			\end{align*}
			which satisfies the initial value problem
			\begin{align*}
				\partial_t\pi_i(t)=-v_i(\pi_i(t)), \qquad \pi_i(0)=x_0.
			\end{align*}
			An argument similar to that in \cite{Tamas} shows that the operator 
			\begin{eqnarray*}
				A:=A_m, \; D(A):= \{f\in W^{1,p}(0,\ell;\R^n):\; f(0)=0\},
			\end{eqnarray*}
			generates a C$_0$-semigroup $\T:=(T(t))_{t\ge 0} $ on $X_p$ given by
			\begin{eqnarray*}%\label{Semigoup-of-A}
				(T(t)f)_{i}(x)=
				\begin{cases}
					e^{\xi_i(\pi_i(t))}f_{i}(\pi_i(t)), & {\rm if } \; t \leq \tau_{i}(x),
					\\
					0,& {\rm if\; not },
				\end{cases}
			\end{eqnarray*}
			for all $i\in \{1,\ldots,n\}$, $t\ge 0$, $f\in X_p$, and a.e. $x\in [0,\ell]$. Clearly, $\T$ is positive and $\rho(A)=\C$. Moreover, the operator $G$ defined in \eqref{trace-operators} is clearlry onto, which ensures the existence of the Dirichlet operator $\D_\la$ associated with $(A_m,G)$. A direct computation then yields that
			\begin{align}\label{Dirichlet-D}
				(\D_\la d)(x) ={\rm diag}\left(e^{\xi_i(x)-\la \tau_i(x)}\right)_{i=1,\ldots,n}d,
			\end{align}
			for all $ d\in \R^n$, $\la \in \C$, and a.e. $x\in [0,\ell]$. We define
			\begin{align*}
				\B:=-A_{-1}\D_0: \R^n\to (X_p)_{-1}^{A}.
			\end{align*}
			Using Laplace transform argument, one obtains that the input-maps of $(A,\B)$ are given by 
			\begin{eqnarray*}%\label{control-maps-B}
				\left(\Phi^{A,\B}_t g\right)_i(x)=\begin{cases} e^{\xi_{i}(x)} g_{i}(t-\tau_{i}(x)),& {\rm if }\; t\ge \tau_i(x),\cr 0,&{\rm if\, not,}\end{cases}
			\end{eqnarray*}
			for all $i\in \{1,\ldots,n\}$, $t\ge 0$, $g\in L^p(\R_+,\R^n)$, and a.e. $x\in [0,\ell]$. By (c) and the change of variables $\sigma=t-\t_i(x)$, we obtain
			\begin{align*}
				\left\Vert \Phi^{A,\B}_tg\right\Vert^p_{X_p}&\le \sum_{i=1}^{n}\int_{0}^{\min\{\ell,\t_i^{-1}(t)\} }e^{p\xi_i(x)}g_i^p(t-\tau_i(x)) dx\\
				&\le v_* e^{p\ell \tfrac{k_*}{\nu}}\Vert g\Vert_{L^p([0,t];\R^n)}^p,
			\end{align*}
			for all $t\ge 0$ and $g\in L^p_+(\R_+,\R^n)$, where $k_*:=\sup_{1\le i\le n}\Vert k_{i}\Vert_\infty$ and $v_*:=\sup_{1\le i\le n}\Vert v_i\Vert_\infty$. Thus, $(A,\B)$ is positive $L^p$-admissible. Moreover, we set $C:=P\vert_{D(A)}$. Then, for $ \max_{i=1,\ldots,n}\tau_i(\ell)\le \alpha$ and $f\in D_+(A)$, we have 
			\begin{align*}
				\int_{0}^{\alpha} \Vert C T(t)f\Vert^p dt &=\sum_{i=1}^{n} \int_{0}^{\tau_i(\ell)} \left(e^{\xi_i(\pi_i(t))}f_{i}(\pi_i(t))\right)^p dt\\
				%	&\le \sum_{i=1}^{n} \int_{0}^{1}  f_i^p(x) dx \\
				&\le  \frac{1}{\nu} e^{p\ell \tfrac{k_*}{\nu}} \Vert f\Vert_{X_p}^p,
			\end{align*}
			where we used the change of variables $x=\pi_i(t)$ and the hypothesis (b) to obtain the above estimation. Thus, $(C,A)$ is positive $L^p$-admissible.
			\vspace{.1cm}
			
			Now, for $\alpha\ge \max_{i=1,\ldots,n}\tau_i(\ell)$, we define the operator:
			\begin{eqnarray*}
				(\mathbb{F}^{A,\B,C}g)(t):=P\Phi^{A,\B}_tg,\quad t\in[0,\alpha],g\in W^{1,p}_{0,+}([0,\alpha];\R^n).
			\end{eqnarray*}
			Then, for any $i\in \{1,\ldots,n\}$ and any $g\in W^{1,p}_{0,+}([0,\alpha];\R^n)$,
			\begin{eqnarray*}%\label{input-output-A}
				(\mathbb{F}^{A,\B,C}g)_i(t)= \begin{cases} e^{\xi_{i}(\ell)} g_{i}(t-\tau_{i}(\ell)),& {\rm if }\; t\ge \tau_i(\ell),\cr 0,&{\rm if\, not.}\end{cases}
			\end{eqnarray*}
			Using the H\"{o}lder inequality, we obtain
			\begin{eqnarray*}
				%\begin{array}{lll}
				\displaystyle\int_{0}^{\alpha} \Vert (\mathbb{F}^{A,\B,C}g)(t)\Vert^p dt&= \displaystyle\sum_{i=1}^{n}\int_{\tau_i(\ell)}^{\alpha} e^{p\xi_i(\ell)}g_i^p(t-\tau_i(\ell)) dt\\&\le  e^{p\ell\tfrac{k_*}{\nu}} \Vert g\Vert_{L^p([0,\alpha];\R^n)}^p,
				%	\end{array}
		\end{eqnarray*}
	for all $g\in W^{1,p}_{0,+}([0,\alpha];\R^n)$. Thus,
\begin{align*}%\label{F-constant}
\Vert\mathbb{F}^{A,\B,C}\Vert_{\calL(L^p(0,\alpha;\R^n))}\le  e^{\ell\tfrac{k_*}{\nu}}.
\end{align*} 
Therefore, $(A,\B,C)$ is a positive $L^p$-well-posed triple on $\R^n,X_p,\R^n$. 

Next, let $s,\la\in \C$ such that $s\ne \la$. Using \eqref{Dirichlet-D} and the resolvent equation, we obtain
\begin{align*}
\lim_{\la \to +\infty} \la CR(\la,A)\D_s d=P\D_sd, \qquad\forall d\in \R^n.
\end{align*}
Thus, Range$(\D_s)\subset D(C_\Lambda)$ and $(C_\Lambda)\vert_{D(A_m)}=P$. Therefore, $\mathbf{H}(\la)=C_\Lambda \D_\la=P\D_\la$. Using \eqref{trace-operators} and \eqref{Dirichlet-D}, one obtains that
\begin{align*}%\label{PDlambda}
\lim_{\la \to +\infty} \Vert P\D_\la\Vert_{\calL(\R^n)}=0,
\end{align*} 
and hence $(A,\B,C)$ is a positive $L^p$-well-posed uniformly regular triple with zero feedthrough.

All assumptions of Corollary \ref{Cor.11} hold for the triple $(A,\B,C)$, and by its claim, the time-delayed linear hyperbolic system $(\Sigma_{\mathsf{Hyp}})$ is well-posed. Furthermore, according to Corollary \ref{Cor1}, $(\Sigma_{\mathsf{Hyp}})$ is uniformly exponentially stable if and only if $r(\Xi  L(\varepsilon_0))<1$, where $\Xi:={\rm diag}\left(e^{\xi_i(\ell)}\right)_{i=1,\ldots,n}$. 
\vspace{.1cm}

Let us now illustrate our results by two specific examples of $\eta$:
\begin{example}
Consider the case where $\eta$ is given by:
\begin{align*}
	\eta(\theta)=\begin{cases}
		\Lambda,& {\rm if}\;\theta=-r,\\
		0,& {\rm if}\; \theta\in (-r,0],
	\end{cases}
\end{align*}  %the entries of the matrix $K$ are negative
where $\Lambda\in \R^{n,n}_+$. By Corollary \ref{Cor3} and Remark \ref{unifrom-stability}, system $(\Sigma_{\mathsf{Hyp}})$ is uniformly exponentially stable if and only if $r(\Xi\Lambda)<1$. Furthermore,
\begin{enumerate}
	\item If $k_{i}(x)\le 0$ for all $i\in \{1,\ldots,n\}$ and a.e. $x\in [0,\ell]$, then $r(\Xi\Lambda)\le r(
	\Lambda)$. Therefore, for $r(\Lambda)<1$, we obtain that $(\Sigma_{\mathsf{Hyp}})$ is uniformly exponentially stable.
	\item If $k_{i}\equiv 0$ for all $i\in \{1,\ldots,n\}$, then $(\Sigma_{\mathsf{Hyp}})$ is uniformly exponentially stable if and only if $r(\Lambda)<1$.
\end{enumerate}
\end{example}

\begin{example}
Consider the case where $\eta$ is given by:
\begin{align*}
	\eta:=\Lambda\otimes \xi.
\end{align*}  
Here, $\Lambda\in \R^{n,n}_+$ and $\xi(\theta):=\int_{-r}^{\theta}e^{\vartheta s}ds$ for all $\theta\in [-r,0]$, where $\vartheta\in\R$ is a parameter. Since $\xi$ is an increasing function and $\Lambda$ is a real positive matrix, then $\eta$ is also an increasing operator-valued function. Moreover, we have ${\rm Var}(\eta;-r,0)\le \Vert\Lambda\Vert \xi(0)$. By Corollary \ref{Cor2} and Remark \ref{unifrom-stability}, the time-delayed hyperbolic system $(\Sigma_{\mathsf{Hyp}})$ is uniformly exponentially stable if and only if $r(\Xi\eta(0))<1$. If all entries of the matrix $K$ are negative, then the spectral radius $r(\Xi\eta(0))$ can be estimated as follows
\begin{eqnarray*}
	\begin{cases}
		r(\Xi\eta(0))<1,& {\rm if} \quad \vartheta<0 \quad {\rm and}\quad 1-e^{-r\vartheta}>\frac{\vartheta}{\Vert\Lambda\Vert},\\
		r(\Xi\eta(0))<1,& {\rm if} \quad \vartheta>0 \quad {\rm and}\quad 1-e^{-r\vartheta}<\frac{\vartheta}{\Vert\Lambda\Vert},\\
		r(\Xi\eta(0))<1,& {\rm if} \quad \vartheta=0 \quad {\rm and}\quad \Vert\Lambda\Vert<\frac{1}{r}.
	\end{cases}
\end{eqnarray*}
\end{example}

\begin{remark}
For the case $p=2$, Corollary \ref{C.1.A} can significantly simplify the well-posedness analysis of system $(\Sigma_{\mathsf{Hyp}})$. This is because both $(A,\B,C)$ and $(Q,\Theta,\L)$ are positive $L^2$-well-posed uniformly regular triples with zero feedthrough. Corollary \ref{C.1.A} can also simplify the proof of positive semigroup generation in \cite[Prop. 5.3 and Thm. 5.4]{BBH}.
\end{remark}

\section{Conclusion}
In this paper, we investigated the well-posedness and stability of a novel class of LTI delay differential equations called \emph{boundary delay equations}. In general, the stability of an abstract delay system cannot be determined by its spectral bound. To tackle this problem, we took advantage of the properties of positive C$_0$-semigroups on Banach lattices. In particular, we showed that the positivity assumption allows us to determine the exponential stability of the boundary delay equations by checking whether the spectral bound is negative. This assumption also enables a relaxation of the conditions for the delay-free initial/boundary-value problem, thereby generalizing well-known criteria for the generation of domain perturbations of semigroup generators \cite{Barbieri,BBH,HMR,MR4635044,El,El1,PiMa}. Furthermore, these findings enable us  to extend and refine existing exponential stability results for positive hyperbolic systems presented in \cite{BBH,Zhang}.

\appendix
\section{Appendix}\label{Appendix}

In this appendix, we present generation and stability theorems for abstract domain perturbations of positive semigroup generators. To explain the general setup, we consider the operator $\calA:D(\calA)\subseteq X\to X$ defined by 
\begin{eqnarray}\label{calA}
\calA:=A_m,\qquad D(\calA):=\left\{x\in D(A_m):\; G x=P x\right\},
\end{eqnarray}
where
\begin{itemize}
\item  $A_m:D(A_m)\subseteq X\to  X$ is a closed and densely defined ``\emph{differential}'' operator;
\item $G,P \in \calL(D(A_m),V)$ are \emph{boundary} operators;
\item $X$ and $V$ are Banach spaces, referred to as the \emph{state} and \emph{boundary} spaces, respectively.
\end{itemize}
Below, we establish sufficient conditions on $A_m$, $G$, and $P$ that guarantee $\mathcal{A}$ generates an exponential stable, positive C$_0$-semigroup on $X$. Throughout, $X$ and $V$ are Banach lattices. We recall the following result (see \cite[Rem. 2.2]{BJVW} and \cite[Prop. 4.3]{SGPS}).
\begin{lemma}
Suppose Assumption \ref{Assp11} holds with $\T$ positive. Then, the following statements are equivalent:
\begin{itemize}
	\item[(i)] $\B$ is positive.
	\item[(i)] $\D_\la$ is positive for all $\la>s(A)$.
	\item[(ii)] $\D_\la$ is positive for all sufficiently large $\la\in \R$.
\end{itemize}
\end{lemma}

We state the following perturbation theorem, which provides mild sufficient conditions under which $\mathcal{A}$ generates a positive semigroup.
\begin{theorem}\label{theorem-ABC}
Let Assumption \ref{Assp11} be satisfied. Suppose that $\T$ is positive, $P$ is positive, and $\D_\la $ is positive for all sufficiently large $\la\in \R$. Assume that 
\begin{itemize}
	\item[(i)] $(A,\B,\P)$ is a positive $L^p$-well-posed weakly regular triple with zero feedthrough, where $\P:=P\vert_{D(A)}$; 
	\item[(ii)] $I_{V}$ is an admissible positive feedback operator for $(A,\B,\P)$.
\end{itemize}
Then, the operator $\calA$, defined by \eqref{calA}, generates a positive C$_0$-semigroup $\calT:=(\calT(t))_{t\ge 0}$ on $ X $ given by $\calT(t)=\T^{cl}(t)$ for all $t\ge 0$.
\end{theorem}
\begin{proof}
The proof follows from \cite[Thm. 4.1]{HMR} and Theorem \ref{Thoerem}.
\end{proof}

\begin{remark}
Theorem \ref{theorem-ABC} is a slight generalization of \cite[Thm. 3.2]{BBH}. The 
positivity of the triple $(A,\B,\P)$ allows us to weaken its regularity. In addition, the assumption that $r(P\D_{\la_0})<1$ for some $\la_0>s(A)$  is no longer needed to prove the positivity of the semigroup $\calT$. 
\end{remark}

\begin{corollary}\label{C.1.A}
Consider the setting of Theorem \ref{theorem-ABC} with $X$ and $V$ are Hilbert lattices. If $(A,\B,\P)$ is a positive $L^2$-well-posed uniformly regular triple with zero feedthrough, then the operator $\calA$ generates a positive C$_0$-semigroup $\calT:=(\calT(t))_{t\ge 0}$ on $ X $.
\end{corollary}
\begin{proof}
According to Definition \ref{regurality}, we have 
\begin{align*}
	\lim_{\la \to +\infty} \Vert P\D_\la\Vert_{\calL(V)}=0,
\end{align*}
since $\mathbf{H}(\la)=P\D_\la$ for all $\la>\omega_0(A)$. Thus, there exists $\la_0\in \R$ sufficiently large such that $r(P\D_{\la_0})<1$. From Remark \ref{Remark-feedabck-admissibility}, we conclude that the identity operator $I_{V}$ is an admissible positive feedback operator for $(A,\B,\P)$. The claim then follows from Theorem \ref{theorem-ABC}.    
\end{proof}

We now recall important stability properties of semigroups (see, e.g., \cite[Sect. 3]{KN}).
\begin{definition}\label{Def.A}
A semigroup $\T$ is \emph{strongly exponential stable} if there exists $\alpha>0$ such that
\begin{align*}
	\Vert T(t)x\Vert_{X} \leq M_{x}e^{-\alpha t}\Vert x\Vert_{X}, 
\end{align*}
for all $x\in D(A)$, $ t\ge 0$, and a constant $M_x>0$. If $M_x$ can be chosen independently of $x$, then $\T$ is called \emph{uniformly exponentially stable}.
\end{definition}
We note that uniform exponential stability of $\T$ is equivalent to $\omega_0(A)<0$. Moreover, when $\T$ is positive, \emph{strong exponential stability} holds if $s(A)<0$; see, e.g., \cite[Thm. A]{KN}.
\vspace{.1cm}

The following theorem provides a spectral characterization of the exponential stability for the semigroup generated by $\calA$.
\begin{theorem}\label{stability-criterion}
Under the assumptions of Theorem \ref{theorem-ABC}, the following equivalence holds
\begin{align*}
	s(\calA)<0  \iff s(A)<0 \quad {\rm and} \quad r(P\D_0)<1.
\end{align*}
In particular, if either of these equivalent conditions is satisfied, the semigroup $\calT$ is strongly exponentially stable.
\end{theorem}
\begin{proof}
It follows from Theorems \ref{Thoerem} and \ref{theorem-ABC} that 
\begin{align}\label{Eq (A.1)}
	s(\calA)\ge s(A).
\end{align}
If $s(\calA)<0$, then by \eqref{Eq (A.1)} one obtains that $s(A)<0$, and therefore the results follows from \cite[Thm. 3.4]{BBH}.
\end{proof}	
				
				\section*{}
				{\bf Acknowledgments}
				\vspace{.1cm}
				
				The authors thank the Associate Editor and the anonymous reviewers for their careful review and insightful comments, which have greatly improved this paper.
				\bibliographystyle{plain}
				\bibliography{cas-refs3}

@article{AABM,
  author  = {Ait Benhassi, E.M. and Ammari, K. and Boulite, S. and Maniar, L.},
  title   = {Feedback stabilization of a class of evolution equations with delay},
  journal = {J. Evol. Equ.},
  volume  = {9},
  pages   = {103--121},
  year    = {2009}
}

@book{CHARALAMBOS,
  author    = {Aliprantis, C.D. and Burkinshaw, O.},
  title     = {Positive Operators},
  publisher = {Springer},
  address   = {Dordrecht, Netherlands},
  year      = {2006},
  doi       = {https://doi.org/10.1007/978-1-4020-5008-4}
}

@article{Arendt,
  author  = {Arendt, W.},
  title   = {Resolvent positive operators},
  journal = {Proc. London Math. Soc.},
  volume  = {54},
  pages   = {321--349},
  year    = {1987}
}

@article{SGPS,
author = {Arora, S. and Gl{\"u}ck, J. and Paunonen, L. and Schwenninger, F.L.},
title = {Limit-case admissibility for positive infinite-dimensional systems},
journal = {J. Differ. Equ.},
volume = {440},
pages = {113435},
year = {2025},
issn = {0022-0396},
doi = {https://doi.org/10.1016/j.jde.2025.113435}
}

@misc{SGF,
  author       = {Arora, S. and Gl{\"u}ck, J. and Schwenninger, F.},
  title        = {The lattice structure of negative {S}obolev and extrapolation spaces},
  howpublished = {arXiv e-print},
  year         = {2024},
  eprint       = {2404.02116},
  eprinttype   = {arxiv},
  doi          = {https://doi.org/10.48550/arXiv.2404.02116}
}

@article{Barbieri,
  title = {Perturbations of positive semigroups factorized via {AM}- and {AL}-spaces},
  author = {Barbieri, A. and Engel, K.-J.},
  journal = {J. Evol. Equ.},
  volume = {25},
  number = {25},
  year = {2025},
  doi = {https://doi.org/10.1007/s00028-024-01049-3}
}

@book{BaCo,
  author    = {Bastin, G. and Coron, J.M.},
  title     = {Stability and boundary stabilization of 1-D hyperbolic systems},
  series    = {Progress in Nonlinear Differential Equations and Their Applications},
  publisher = {Springer},
  address   = {New York, NY, USA},
  year      = {2016}
}

@article{Boujijane2024,
  author  = {Boujijane, S. and Boulite, S. and Halloumi, M. and Maniar, L. and Rhandi, A.},
  title   = {Well-posedness and asynchronous exponential growth of an age-weighted structured fish population model with diffusion in {$L^1$}},
  journal = {J. Evol. Equ.},
  year    = {2024},
  volume  = {24},
  number  = {14},
  doi     = {https://doi.org/10.1007/s00028-023-00942-7}
}

@article{BP1,
  author  = {B{\'a}tkai, A. and Piazzera, S.},
  title   = {Semigroups and linear partial differential equations with delay},
  journal = {J. Math. Anal. Appl.},
  volume  = {264},
  pages   = {1--20},
  year    = {2001}
}

@book{BP2,
  author    = {B{\'a}tkai, A. and Piazzera, S.},
  title     = {Semigroups for Delay Equations},
  series    = {Research Notes in Mathematics},
  volume    = {10},
  publisher = {AK Peters, Ltd.},
  address   = {Wellesley, MA},
  year      = {2005}
}

@article{BJVW,
  author  = {B{\'a}tkai, A. and Jacob, B. and Voigt, J. and Wintermayr, J.},
  title   = {Perturbation of positive semigroups on {AM}-spaces},
  journal = {Semigroup Forum},
  volume  = {96},
  pages   = {33--47},
  year    = {2018}
}

@misc{BFR1,
  author    = {B{\'a}tkai, A. and Fijav{\v z}, M.K. and Rhandi, A.},
  title     = {Abstract boundary delay systems and application to network flow},
 howpublished = {arXiv e-print},
  year         = {2025},
  eprint       = {2503.08809},
  eprinttype   = {arxiv},
  doi          = {https://doi.org/10.48550/arXiv.2503.08809}
}

@book{BFR,
  author    = {B{\'a}tkai, A. and Fijav{\v z}, M.K. and Rhandi, A.},
  title     = {Positive Operator Semigroups: From Finite to Infinite Dimensions},
  publisher = {Birkh{\"a}user-Verlag},
  address   = {Basel},
  year      = {2016}
}

@article{BBH,
  author  = {Boulouz, A. and Bounit, H. and Hadd, S.},
  title   = {Feedback theory approach to positivity and stability of evolution equations},
  journal = {Syst. Control Lett.},
  volume  = {161},
  pages   = {105167},
  year    = {2022}
}

@book{CZ,
  author    = {Curtain, R. and Zwart, H.},
  title     = {An Introduction to Infinite-Dimensional Linear Systems Theory},
  volume    = {21},
  series    = {Texts in Applied Mathematics},
  publisher = {Springer-Verlag},
  address   = {New York},
  year      = {1995}
}

@article{Datko,
  author  = {Datko, R.},
  title   = {Not all feedback stabilized hyperbolic systems are robust with respect to small time delay in their feedbacks},
  journal = {SIAM J. Control Optim.},
  volume  = {26},
  pages   = {697--713},
  year    = {1988}
}

@article{DLP,
  author  = {Datko, R. and Lagnese, J. and Polis, M.P.},
  title   = {An example on the effect of time delays in boundary feedback stabilization of wave equations},
  journal = {SIAM J. Control Optim.},
  volume  = {24},
  pages   = {152--156},
  year    = {1986}
}

@article {MR4635044,
    AUTHOR = {El Gantouh, Yassine},
     TITLE = {Boundary approximate controllability under positivity
              constraints of infinite-dimensional control systems},
   JOURNAL = {J. Optim. Theory Appl.},
    VOLUME = {198},
      YEAR = {2023},
     PAGES = {449--478},
}

@misc{El,
  author       = {El Gantouh, Y.},
  title        = {Positivity of infinite dimensional linear systems},
  howpublished = {arXiv e-print},
  year         = {2023},
  eprint       = {2208.10617},
  eprinttype   = {arxiv},
  doi         = {https://doi.org/10.48550/arXiv.2208.10617}
}

@article{El1,
  author  = {El Gantouh, Y.},
  title   = {Well-posedness and stability of a class of linear systems},
  journal = {Positivity},
  volume  = {28},
  number  = {16},
  year    = {2024},
  doi     = {https://doi.org/10.1007/s11117-024-01035-6}
}

@book{EN,
  author    = {Engel, K.J. and Nagel, R.},
  title     = {One-Parameter Semigroups for Linear Evolution Equations},
  publisher = {Springer-Verlag},
  address   = {New York},
  year      = {2000}
}

@article{Gr,
  author  = {Greiner, G.},
  title   = {Perturbing the boundary conditions of a generator},
  journal = {Houston J. Math.},
  volume  = {13},
  pages   = {213--229},
  year    = {1987}
}

@article{PiMa,
  author  = {Gwi{\'z}d{\'z}, P. and Tyran-Kami{\'n}ska, M.},
  title   = {Positive semigroups and perturbation of boundary conditions},
  journal = {Positivity},
  volume  = {23},
  pages   = {921--939},
  year    = {2019}
}

@article{HIR,
  author  = {Hadd, S. and Idrissi, A. and Rhandi, A.},
  title   = {The regular linear systems associated to the shift semigroups and application to control delay systems},
  journal = {Math. Control Signals Syst.},
  volume  = {18},
  pages   = {272--291},
  year    = {2006}
}

@article{HMR,
  author  = {Hadd, S. and Manzo, R. and Rhandi, A.},
  title   = {Unbounded perturbations of the generator domain},
  journal = {Discrete Contin. Dyn. Syst.},
  volume  = {35},
  pages   = {703--723},
  year    = {2015}
}

@book{Hale,
  author    = {Hale, J.K.},
  title     = {Functional Differential Equations},
  publisher = {Springer-Verlag},
  address   = {New York},
  year      = {1971}
}

@article{HK,
  author  = {Hashimoto, T. and Krstic, M.},
  title   = {Stabilization of reaction-diffusion equations with state delay using boundary control input},
  journal = {IEEE Trans. Autom. Control},
  volume  = {61},
  pages   = {4041-4047},
  year    = {2016}
}

@article{KN,
  author  = {Kerscher, W. and Nagel, R.},
  title   = {Asymptotic behavior of one-parameter semigroups of positive operators},
  journal = {Acta Appl. Math.},
  volume  = {2},
  pages   = {297--309},
  year    = {1984}
}

@article{JGH,
  author  = {Jiang, W. and Guo, F. and Huang, F.},
  title   = {Well-posedness of linear partial differential equations with unbounded delay operators},
  journal = {J. Math. Anal. Appl.},
  volume  = {293},
  pages   = {310--328},
  year    = {2004}
}

@article{Tamas,
  author  = {M{\'a}trai, T. and Sikolya, E.},
  title   = {Asymptotic behavior of flows in networks},
  journal = {Forum Math.},
  volume  = {19},
  pages   = {429--461},
  year    = {2007}
}

@incollection{MVo,
  author    = {Maniar, L. and Voigt, J.},
  title     = {Linear delay equations in the {$L^p$}-context},
  booktitle = {Recent Contributions to Evolution Equations},
  publisher = {Marcel Dekker},
  address   = {New York},
  pages     = {319--330},
  year      = {2003},
  series    = {Lecture Notes in Pure and Applied Mathematics}
}

@book{Nagel,
  editor    = {Nagel, R.},
  title     = {One-Parameter Semigroups of Positive Operators},
  series    = {Lecture Notes in Mathematics},
  publisher = {Springer},
  address   = {Berlin},
  year      = {1986}
}

@article{Salam,
  author  = {Salamon, D.},
  title   = {Infinite-dimensional linear system with unbounded control and observation: a functional analytic approach},
  journal = {Trans. Amer. Math. Soc.},
  volume  = {300},
  pages   = {383--431},
  year    = {1987}
}

@book{Schaf,
  author    = {Schaefer, H.H.},
  title     = {Banach Lattices and Positive Operators},
  publisher = {Springer-Verlag},
  address   = {Berlin-Heidelberg},
  year      = {1974}
}

@book{Staf,
  author    = {Staffans, O.J.},
  title     = {Well-posed Linear Systems},
  publisher = {Cambridge University Press},
  address   = {Cambridge},
  year      = {2005}
}

@article{Webb,
  author  = {Webb, G.},
  title   = {Functional differential equations and nonlinear semigroups in {$L^p$}-spaces},
  journal = {J. Differ. Equ.},
  volume  = {29},
  pages   = {71--89},
  year    = {1976}
}

@article{WC,
  author  = {Weiss, G.},
  title   = {Admissibility of unbounded control operators},
  journal = {SIAM J. Control Optim.},
  volume  = {27},
  pages   = {527--545},
  year    = {1989}
}

@book{TW,
  author    = {Tucsnak, M. and Weiss, G.},
  title     = {Observation and Control for Operator Semigroups},
  publisher = {Birkh{\"a}user},
  address   = {Basel, Boston, Berlin},
  year      = {2009}
}

@article{WF,
  author  = {Weiss, G.},
  title   = {Regular linear systems with feedback},
  journal = {Math. Control Signals Syst.},
  volume  = {7},
  pages   = {23--57},
  year    = {1994}
}

@article{WR,
  author  = {Weiss, G.},
  title   = {Transfer functions of regular linear systems. {Part I}: Characterization of regularity},
  journal = {Trans. Amer. Math. Soc.},
  volume  = {342},
  pages   = {827--854},
  year    = {1994}
}

@article{Weis,
  author  = {Weis, L.},
  title   = {The stability of positive semigroups on {$L^p$}-spaces},
  journal = {Proc. Amer. Math. Soc.},
  volume  = {123},
  pages   = {3089--3094},
  year    = {1995}
}

@book{Wu,
  author    = {Wu, J.},
  title     = {Theory and Applications of Partial Functional Differential Equations},
  publisher = {Springer-Verlag},
  address   = {New York},
  series    = {Applied Mathematical Sciences},
  volume    = {119},
  year      = {1996}
}

@article{Xu,
  author  = {Xu, G.Q.},
  title   = {Resolvent family for the evolution process with memory},
  journal = {Math. Nachr.},
  volume  = {296},
  pages   = {2626--2656},
  year    = {2023}
}

@article{Zhang,
  author  = {Zhang, L. and Prieur, C.},
  title   = {Necessary and sufficient conditions on the exponential stability of positive hyperbolic systems},
  journal = {IEEE Trans. Autom. Control},
  volume  = {62},
  pages   = {3610--3617},
  year    = {2017}
}

@article{StaffansWeiss2002,
	author = {Staffans, O. and Weiss, G.},
	title = {Transfer functions of regular linear systems. {Part II}: The system operator and the {Lax-Phillips} semigroup},
	journal = {Trans. Amer. Math. Soc.},
	volume = {354},
	pages = {3229--3262},
	year = {2002}
}

@book{fridman2014introduction,
  author    = {Fridman, Emilia},
  title     = {Introduction to Time-Delay Systems: Analysis and Control},
  publisher = {Birkh\"{a}user},
  year      = {2014},
  series    = {Systems \& Control: Foundations \& Applications},
  address   = {Cham, Switzerland},
  isbn      = {978-3-319-09392-5},
  doi       = {10.1007/978-3-319-09393-2},
  pages     = {364}
}

				%  \bibliographystyle{elsarticle-num}
				%\bibliography{dada}
			\end{document}